%% file: arXiv-SACOBRA15.tex
\PassOptionsToPackage{table}{xcolor}
\pdfoutput=1				% for arXiv

\let\oldvec\vec
\documentclass{elsarticle}
\let\vec\oldvec
\usepackage[latin1]{inputenc}
\usepackage[english]{babel}
\usepackage{xspace}			% xspace
\usepackage{paralist}		% compactitem usw.
\usepackage{longtable} %long tables multiple pages
\usepackage{multirow}
\usepackage{xcolor}% vielleicht versteckt in einem anderen Paket

\usepackage{eurosym}
\usepackage{amsmath}		% bmatrix and the like
\DeclareMathOperator*{\med}{median}

\usepackage{lipsum}
\usepackage{url}

\makeatletter
\usepackage{subfigure}
\usepackage{booktabs}						% toprule, midrule, bottomrule
\usepackage{color}
\usepackage{algorithm}
\usepackage{algpseudocode}
\usepackage{amsfonts}
\usepackage{amssymb} 						% \mathbb
\usepackage{algpseudocode}
\usepackage[bookmarks=false]{hyperref}						% \href
\usepackage{stackrel}
\usepackage{enumerate}					% roman enumeration (i)
\usepackage{tikz}
\usetikzlibrary{patterns}
\usetikzlibrary{shapes,arrows} %flowchart
\usetikzlibrary{positioning}
\usetikzlibrary{arrows}
\usetikzlibrary{shapes.geometric, arrows,positioning}
\DeclareGraphicsExtensions{pdf,png}

\newcommand{\SB}[1]{\textcolor{magenta}{\small{(\textbf{SB:} \textit{#1})}}}
\newcommand{\addOK}[1]{#1}   % just a helper to turn an \addOK into an unconditionally accepted add 
\newcommand{\deletedOK}[1]{} 
\newcommand{\replacedOK}[2]{#1}

\newcommand{\mat}[1]{\mathbf{#1}}
\makeatletter
\tikzset{nomorepostaction/.code=\let\tikz@postactions\pgfutil@empty}
\makeatother
\journal{Applied Soft Computing}

\begin{document}

\pagestyle{headings}  % switches on printing of running heads

%% Title, authors and addresses

%% use the tnoteref command within \title for footnotes;
%% use the tnotetext command for theassociated footnote;
%% use the fnref command within \author or \address for footnotes;
%% use the fntext command for theassociated footnote;
%% use the corref command within \author for corresponding author footnotes;
%% use the cortext command for theassociated footnote;
%% use the ead command for the email address,
%% and the form \ead[url] for the home page:
%% \title{Title\tnoteref{label1}}
%% \tnotetext[label1]{}
%% \author{Name\corref{cor1}\fnref{label2}}
%% \ead{email address}
%% \ead[url]{home page}
%% \fntext[label2]{}
%% \cortext[cor1]{}
%% \address{Address\fnref{label3}}
%% \fntext[label3]{}

\begin{frontmatter}

\title{
%Solving the G-problems in less than 500 iterations: Improved efficient constrained optimization by self-adjusting surrogate models  \\
%       OR \\
%Solving the G-problems in less than 500 iterations: Surrogate modeling with adaptive parameter control improves efficient constrained optimization
Solving the G-problems in less than 500 iterations: \\
Improved efficient constrained optimization by surrogate modeling and adaptive parameter control
%			 \\ {\normalsize\WK{decide about title}}
}

%% use optional labels to link authors explicitly to addresses:
%% \author[label1,label2]{}
%% \address[label1]{}
%% \address[label2]{}

\author[CGM]{Samineh~Bagheri}
\author[CGM]{Wolfgang Konen\corref{cor1}}
\ead{\{samineh.bagheri,wolfgang.konen\}@th-koeln.de}
\cortext[cor1]{Corresponding Author}

\address[CGM]{Department of Computer Science, \\ TH K\"{o}ln (University of Applied Sciences), 51643 Gummersbach, Germany}
\author[LEI]{Michael Emmerich} 
\author[LEI]{Thomas B\"ack} 
\ead{\{m.t.m.emmerich,T.H.W.Baeck\}@liacs.leidenuniv.nl}
\address[LEI]{Leiden University, LIACS, \\2333 CA Leiden, The Netherlands}
%divis intelligent solutions GmbH,\\44227 Dortmund\\
%\email{\{baeck,krause,foussette\}@divis-gmbh.de}

%\maketitle

\begin{abstract}
\boldmath
Constrained optimization of high-dimensional numerical problems plays an important role in many scientific and industrial applications. Function evaluations in many industrial applications are severely limited and no analytical information about objective function and constraint functions is available. For such expensive black-box optimization tasks, the constraint optimization algorithm COBRA was proposed, making use of RBF surrogate modeling for both the objective and the constraint functions. COBRA has shown remarkable success in solving reliably complex benchmark problems in less than 500 function evaluations. Unfortunately, COBRA requires careful adjustment of parameters in order to do so. 

In this work we present a new self-adjusting algorithm SACOBRA, which is based on COBRA and capable to achieve high-quality results with very few function evaluations and no parameter tuning. It is shown with the help of performance profiles
%\SB{shall we change this to data profiles?} \WK{I think performance profile is better understandable at this point}
on a set of benchmark problems (G-problems, MOPTA08) that SACOBRA consistently outperforms any COBRA algorithm with fixed parameter setting. We analyze the importance of the several new elements in SACOBRA and find that \replacedOK{each element of SACOBRA plays a role to boost up the overall optimization performance.}{rescaling of input dimensions and random restart are the most important elements, followed by a self-adjusting output transform mechanism.} 
We discuss the reasons behind and get in this way a better understanding of high-quality RBF surrogate modeling. 
\end{abstract}

\begin{keyword}
%% keywords here, in the form: keyword \sep keyword

%% PACS codes here, in the form: \PACS code \sep code

%% MSC codes here, in the form: \MSC code \sep code
%% or \MSC[2008] code \sep code (2000 is the default)

optimization; constrained optimization; expensive black-box optimization; radial basis function; self-adjustment

\end{keyword}

\end{frontmatter}

\pagebreak

%%%%%%%%%%%%%%%%%%%%%%%%%%%%%%%%%%%%%%%%%%%%%%%%%%%%%%%%%%%%%%%%%%%%%%%%%%%%%%%%%%%%%%%%%%%%%%%%%%%%
\section{Introduction}
%%%%%%%%%%%%%%%%%%%%%%%%%%%%%%%%%%%%%%%%%%%%%%%%%%%%%%%%%%%%%%%%%%%%%%%%%%%%%%%%%%%%%%%%%%%%%%%%%%%%
\label{Sec:Introduction}

Real-world optimization problems are often subject to constraints, restricting the feasible region to a smaller subset of the search space. It is the goal of constraint optimizers to avoid infeasible solutions and to stay in the feasible region, in order to converge to the optimum. However, the search in constraint black-box optimization can be difficult, since \replacedOK{we usually have no a-priori knowledge about the feasible region and the fitness landscape}{knowledge about the size of the feasible region and the fitness landscape is usually unknown}.
This problem even turns out to be harder, when only a limited number of function evaluations is allowed for the search. However, in industry good solutions are requested in very restricted time frames. An example is the well-known benchmark MOPTA08~\cite{jones2008large}.

In the past different strategies have been proposed to handle constraints. E.~g., repair methods try to guide infeasible solutions into the feasible area. Penalty functions give a negative bias to the objective function value, when constraints are violated. Many constraint handling methods are available in the scientific literature, but often demand for a large number of function evaluations~(e.~g., results in \cite{runarsson2000stochastic,kramer2006three}). %This can be especially problematic for real-world optimization tasks, which sometimes incorporate expensive simulations. 
Up to now, only little work has been devoted to \emph{efficient} constraint optimization (severely reduced number of function evaluations). A possible solution in that regard is to use \emph{surrogate models} for the objective and the constraint functions. While the real function might be expensive to evaluate, evaluations on the surrogate functions are usually cheap. As an example for this approach, the solver COBRA (Constrained Optimization by Radial Basis Function Approximation) was proposed by Regis~\cite{regis2014constrained} and outperforms many other algorithms on a large number of benchmark functions.

In our previous work~\cite{Koch14a,Koch2015a} we have studied a reimplementation of COBRA in \texttt{R}~\cite{Rproject13} enhanced by a new repair mechanism and reported its strengths and weaknesses. Although good results were obtained, each new problem required tedious manual tuning of the many parameters in COBRA. In this paper we follow a more unifying path and present SACOBRA (Self-Adjusting COBRA), an extension of COBRA which starts with the same settings on all problems and adjusts all necessary parameters internally.\footnote{
SACOBRA is available as open-source \texttt{R}-package from CRAN: \href{https://cran.r-project.org/web/packages/SACOBRA}{\url{https://cran.r-project.org/web/packages/SACOBRA}}}
This is \textit{adaptive parameter control} according to the terminology introduced by Eiben et al.~\cite{eiben1999}. We present extensive tests of SACOBRA and other algorithms on a well-known and popular benchmark from the literature: The so-called G-problem or G-function benchmark was introduced by Michalewicz and Schoenauer~\protect\cite{michalewicz1996evolutionary} and Floudas and Pardalos~\cite{Floudas}. It provides a set of constrained optimization problems with a wide range of different conditions.

We define the following research questions for the constrained optimization experiments in this work:
%\SB{I changed the order of research questions. Only if you are also fine with it otherwise I can turn it back.}
%\WK{is o.k.}
\begin{itemize}
	\item[\textbf{(H1)}] Do numerical instabilities occur in RBF surrogates and is it possible to avoid them?
	\item[\textbf{(H2)}] Is it possible with SACOBRA to start with the same initial parameters on all G-problems and to solve them by self-adjusting the parameters on-line?
	\item[\textbf{(H3)}] Is it possible with SACOBRA to solve all G-problems in less than 500 function evaluations?
	%cope with vastly different conditions in constraint optimization (e.~g. G-problems)?
	%\item[(H4)] Gaussian RBF or cubic RBF? Can we find self-adjusting procedure for Gaussian width? \WK{next paper}
\end{itemize}

%%%%%%%%%%%%%%%%%%%%%%%%%%%%%%%%%%%%%%%%%%%%%%%%%%%%%%%%%%%%%%%%%%%%%%%%%%%%%%%%%%%%%%%%%%%%%%%%%%%%
\subsection{Related work}
Following the surveys on constraint optimization given by Michalewicz and Schoenauer~\cite{michalewicz1996evolutionary}, Eiben and Smith~\cite{eiben2003introduction}, Coello Coello~\cite{coello2012constraint}, Jiao et al.~\cite{Jiao2013},  and Kramer~\cite{kramer2010review}, several approaches are available for \textbf{constraint handling}:
\begin{enumerate}[(i)]
	\item unconstrained optimization with a penalty added to the fitness value for infeasible solutions
	\item feasible solution preference methods and stochastic ranking
	\item repair algorithms to resolve constraint violations during the search 
	\item multi-objective optimization, where the constraint functions are defined as additional objectives
\end{enumerate}
A frequently used approach to handle constraints is to incorporate static or dynamic penalty terms (i) in order to stay in the feasible region~\cite{coello2012constraint,kramer2010review,mezura2011constraint}. Penalty functions can be very helpful for solving constrained problems, but their main drawback is that they often require additional parameters for balancing the fitness and penalty terms. Tessema and Yen~\cite{tessema2009} propose an interesting adaptive penalty method which does not need any parameter tuning.

Feasible solution preference methods (ii) \cite{deb2000efficient,mezura2005simple} always prefer feasible solutions to infeasible solutions. They may use too little information from infeasible solutions and risk getting stuck in local optima. Deb~\cite{deb2000efficient} improves this method by introducing a diversity mechanism.
Stochastic ranking~\cite{runarsson2000stochastic,runarsson2005search} is a similar and very successful improvement: With a certain probability %$P_f$
an infeasible solution is ranked not behind, but -- according to its fitness value --  \textit{among} the feasible solutions. 
Stochastic ranking has shown good results on all 11 G-problems. 
%\SB{Not all of the G-functions are introduced by Michalewicz} \WK{o.k.}
However it requires usually a large number of function evaluations (300\,000 and more) and is thus not well suited for efficient optimization.

Repair algorithms (iii) try to transform infeasible solutions into feasible ones~\cite{chootinan2006constraint,zahara2009,Koch2015a}. The work of Chootinan and Chen~\cite{chootinan2006constraint} shows very good results on 11 G-problems, but requires a large number of function evaluations (5\,000 -- 500\,000) as well.

In recent years, multi-objective optimization techniques (iv) have attracted increasing attention for solving constrained optimization problems. The general idea is to treat the constraints as one or more objective functions to be optimized in conjunction with the fitness function. Coello Coello and Montes~\cite{CoelloMontes2002} use Pareto dominance-based tournament selection for a genetic algorithm (GA). Similarly, Venkatraman and Yen~\cite{Venkatra2005} propose a two-phase GA, where the second phase is formulated as a bi-objective optimization problem which uses non-dominated ranking. Jiao et al.~\cite{Jiao2013} use a novel selection strategy based on bi-objective optimization and get improved reliability on a large number of benchmark functions. Emmerich et al.~\cite{emmerich2006} use Kriging models for approximating constraints in a multi-objective optimization scheme.

In the field of \textbf{model-assisted optimization} algorithms for constrained problems, support vector machines~(SVMs) have been used by Poloczek and Kramer~\cite{poloczek2013local}. They make use of SVMs as a classifier for predicting the feasibility of solutions, but achieve only slight improvements. Powell~\cite{powell1994direct} proposes COBYLA, a direct search method which models the objective and the constraints using linear approximation. Recently, Regis~\cite{regis2014constrained} developed COBRA, an efficient solver that makes use of Radial Basis Function~(RBF) interpolation \addOK{to model objective and constraint functions},
%\SB{my suggestion to add this extra explanation}\WK{o.k.}
and outperforms most algorithms in terms of required function evaluations on a large number of benchmark functions.
Tenne and Armfield~\cite{tenne2008} present %nice work where they use 
an adaptive topology RBF network to tackle highly multimodal functions. But they consider only unconstrained optimization.

Most optimization algorithms need their parameter to be set with respect to the specific optimization problem in order to show good performance. Eiben et al.~\cite{eiben1999} introduced a terminology for parameter settings for evolutionary algorithms: They distinguish parameter tuning (before the run) and parameter control (online). Parameter control is further subdivided into predefined control schemes (deterministic),  control with feedback from the optimization run (adaptive), or  control where the parameters are part of the evolvable chromosome (self-adaptive). 

Several papers deal with \textbf{adaptive or self-adaptive parameter control} in unconstrained or constrained optimization:
%Work on \textbf{self-adjusting optimization} algorithms:\\
%Cai et al.~\cite{cai2008pso} propose a self-adjusting algorithm for particle swarm optimization. \WK{This work seems not very relevant, only PSO on some unconstrained test functions}  \\
Qin and Suganthan~\cite{qin2005saDE} propose a self-adaptive differential evolution (DE) algorithm.
Brest et al.~\cite{brest2006saDE} propose another self-adaptive DE algorithm. But they do not handle constraints, whereas Zhang et al.~\cite{Zhang2012} describe a constraint-handling mechanism for DE. We will compare later our results with the DE-implementation \texttt{DEoptimR}\footnote{\texttt{R}-package \texttt{DEoptimR}, available from \href{https://cran.r-project.org/web/packages/DEoptimR}{\url{https://cran.r-project.org/web/packages/DEoptimR}}} which is based on both works~\cite{brest2006saDE,Zhang2012}. 
%\WK{Brest et al. give repeat the good definition of deterministic, adaptive and self-adaptive parameter control which goes back to Eiben et al.~\cite{eiben1999}.} \\
Farmani and Wright~\cite{farmani2003saFitness} propose a self-adaptive fitness formulation and test it on all 11 G-problems. They show \replacedOK{comparable results}{ good results similar} to stochastic ranking \cite{runarsson2000stochastic}, but require many function evaluations (above 300\,000) as well.
%\WK{Do we know the number of fe's in \cite{farmani2003saFitness}?}
%\SB{It is described on page 451 that a maximum of 5000 generation is allowed for all the problems and in average after 26 generations first feasible point is appeared. But for G10, 517 generations are required to find the first feasible point and for 3 out of 20 runs no feasible point is found. Also for G05 only 9 out of 20 runs can find a feasible solution (I think in every generation 70 points are evaluated)} 
%\WK{o.k., this translates to a maximum of 5000*70 = 350\,000 fe's}
%
%\WK{Do we know the dimension number used for G02 and G03 in \cite{farmani2003saFitness}?}
%\SB{I guess 20 dimension for G2, because of the optimum value which is  about 0.8. And because she did not discuss anything about G03 I assume she used the commonly used dimension for this problem which is 20.}
%
%\WK{Another remark: Farmani is wrong when saying on p. 451 that RunarssonYao~\cite{runarsson2000stochastic} have in G10 only few feasible runs (6 out of 20). SB found that this is true only for some intermediate results of RunarssonYao, in the end they show an algorithm which produces feasible solutions in all 20 runs.}
%
%Lin~\cite{Lin2013} proposes a novel genetic algorithm based on rough set theory which contains a parameter tuning policy.
Coello Coello~\cite{coello2000saPenalty} and Tessema and Yen~\cite{tessema2009} propose self-adaptive penalty approaches. A survey of self-adaptive penalty approaches is given in \cite{eiben1999}.
%\WK{I did only superfluously read the articles so far, cannot assure quality yet. But Brest et al.\cite{brest2006saDE} seems good quality, was cited in DEoptimR package.}

The area of \textbf{efficient constrained optimization}, that is optimization under severely limited budget of less than 1000 function evaluations, is attracting more and more attention in recent years: 
Regis proposed besides the already mentioned COBRA approach~\cite{regis2014constrained} a trust-region evolutionary algorithm~\cite{regis2015trust} which uses RBF surrogates as well and which exhibits high-quality results on many but not all G-functions in less than 1000 function evaluations.
%\WK{Regis does in \cite{regis2015trust} also transform manually some fitness functions and some constraints. He is further not clear about which of the G-problems are really solved (I assume that G02 is not solved by any of the algorithms shown in the performance profiles.}
%
Jiao et al.~\cite{Jiao2013} propose a self-adaptive selection method to combine informative feasible and infeasible solutions and they formulate it as a  multi-objective problem. Their algorithm can solve some of the G-functions (G08,G11,G12) really fast in less than 500 evaluations, some others are solved with less than 10\,000 evaluations, but the remaining G-functions (G01-G03,G07,G10) require 20\,000 to 120\,000 evaluations to be solved.  
%\WK{see Table 7 in \cite{Jiao2013} on p. 136}
%
Zahara and Kao~\cite{zahara2009} show similar results (1000 -- 20\,000 evaluations) on some G-functions, but they investigate only G04, G08, and G12.
%Why only these 3 G-funcs? No good result for others?
%
To the best of our knowledge there is currently no approach which can solve all 11 G-problems in less than 1\,000 evaluations.
Tenne and Armfield~\cite{tenne2008} present an interesting approach with approximating RBFs to optimize highly multimodal functions in less than 200 evaluations, but their results are only for unconstrained functions and they are not competitive in terms of precision.

The rest of this paper is organized as follows: In Sec.~\ref{Sec:Methods} we present the constrained optimization problem and our methods: the RBF surrogate modeling technique, the COBRA-R algorithm and the SACOBRA algorithm. In Sec.~\ref{Sec:Experiments} we perform a thorough experimental study on analytical test functions and on a real-world benchmark function MOPTA08~\cite{jones2008large} from the automotive domain. We analyze with the help of data profiles the impact of the various SACOBRA elements on the overall performance. The results are discussed in Sec.~\ref{Sec:Discussion} and we give conclusive remarks in Sec.~\ref{Sec:Conclusion}. 

%%%%%%%%%%%%%%%%%%%%%%%%%%%%%%%%%%%%%%%%%%%%%%%%%%%%%%%%%%%%%%%%%%%%%%%%%%%%%%%%%%%%%%%%%%%%%%%%%%%%
\section{Methods}
%%%%%%%%%%%%%%%%%%%%%%%%%%%%%%%%%%%%%%%%%%%%%%%%%%%%%%%%%%%%%%%%%%%%%%%%%%%%%%%%%%%%%%%%%%%%%%%%%%%%
\label{Sec:Methods}

%%%%%%%%%%%%%%%%%%%%%%%%%%%%%%%%%%%%%%%%%%%%%%%%%%%%%%%%%%%%%%%%%%%%%%%%%%%%%%%%%%%%%%%%%%%%%%%%%%%%
\subsection{Constrained optimization}
A constrained optimization problem can be defined by the minimization of an objective function $f$ subject to constraint function(s) $g_1,\ldots,g_m$:
\begin{eqnarray}
 \text{Minimize}    & f(\vec{x}) & \label{eq:minprob}\\
 \text{subject to}  \nonumber \\
										& g_{i}(\vec{x}) \leq 0, \quad &i=1,2,\ldots,m ,  \nonumber \\
										& \vec{x} \in [\vec{a},\vec{b}] &\subset \mathbb{R}^d \nonumber
\end{eqnarray}
In this paper we always consider minimization problems. Maximization problems are transformed to minimization without loss of generality. Problems with equality constraints have them transformed to inequalities first (see Sec.~\ref{Sec:ExpSetup} and Sec.~\ref{sec:limitation}).

%%%%%%%%%%%%%%%%%%%%%%%%%%%%%%%%%%%%%%%%%%%%%%%%%%%%%%%%%%%%%%%%%%%%%%%%%%%%%%%%%%%%%%%%%%%%%%%%%%%%
\subsection{Radial Basis Functions}
\label{Sec:RBF}
The COBRA algorithm incorporates optimization on auxiliary functions, e.~g. regression models over the search space. Although numerous regression models are available, we employ interpolating RBF models~\cite{buhmann2003radial,powell1992theory}, since they outperform other models in terms of efficiency and quality. In this paper we use the same notation as Regis~\cite{regis2014particle}. RBF models require as input a set of design points (a training set):
$n$ points $\vec{u}^{(1)},\ldots,\vec{u}^{(n)} \in \mathbb{R}^d$ are evaluated on the real function $f(\vec{u}^{(1)}),\ldots,f(\vec{u}^{(n)})$. We use an interpolating radial basis function as approximation:
\begin{equation}
%\hat{f}(\vec{x}) = \sum_{i=1}^n \lambda_i \varphi(||\vec{x}^{(+)} - \vec{x}^{(i)}||) + p(\vec{x}), \quad\vec{x} \in \mathbb{R}^d
 s^{(n)}(\vec{x}) = \sum_{i=1}^n \lambda_i \varphi(||\vec{x}       - \vec{u}^{(i)}||) + p(\vec{x}), \quad\vec{x} \in \mathbb{R}^d
\label{eq:RBFfunction}
\end{equation}

Here, $||\cdot||$ is the Euclidean norm, $\lambda_i \in \mathbb{R}$ for $i=1,\ldots,n$, $p(\vec{x}) = c_0 + \vec{c}\,\vec{x}$ is a linear polynomial in $d$ variables with $d+1$ coefficients $\vec{c}\,'=(c_0,\vec{c})^T=(c_0,c_1,\ldots,c_{d})^T \in \mathbb{R}^{d+1}$, and $\varphi$ is of cubic form $\varphi(r) = r^3$. An alternative to cubic RBFs are Gaussian RBFs $\varphi(r) = e^{-r^2/(2\sigma^2)}$ which introduce an additional width parameter $\sigma$.
%Although other choices for $\varphi$ are possible and have been tested in related work, cubic RBF have been shown to be superior in~\cite{wild2013global}. However, as the cubic form showed successful results in earlier studies and has been used by Regis~\cite{regis2014constrained} in the original COBRA algorithm, we decided to integrate it as regression model in our COBRA implementation.

The RBF model fit requires a distance matrix $\mat{\Phi} \in \mathbb{R}^{n\times n}$: $\Phi_{ij} = \varphi(||\vec{u}^{(i)}-\vec{u}^{j}||), i,j=1,\ldots,n$. The RBF model requires the solution of the following linear system of equations:
\begin{equation}
\begin{bmatrix} \mat{\Phi} & \mat{P} \\ \mat{P}^T & \mat{0}_{(d+1)\times(d+1)} \end{bmatrix}
\begin{bmatrix} \vec{\lambda} \\ \vec{c}\,' \end{bmatrix}  = \left[ \begin{array}{c} F \\ 0_{d+1} \end{array} \right]
\label{eq:RBF}
\end{equation}
for the unknowns $\vec{\lambda}, \vec{c}\,'$. Where $\mat{P} \in \mathbb{R}^{n\times(d+1)}$ is a matrix with $(1,\vec{u}^{(i)})$ in its $i$th row, $\mat{0}_{(d+1)\times(d+1)} \in \mathbb{R}^{(d+1)\times(d+1)}$ is a zero matrix, $0_{d+1}$ is a vector of zeros, $F=(f(\vec{u}^{(1)}),	\allowbreak\ldots, \allowbreak f(\vec{u}^{(n)}))^T$, and $\vec{\lambda} = (\lambda_1,\ldots, \lambda_n)^T \in \mathbb{R}^n$. The matrix in Eq.~\eqref{eq:RBF} is invertible if it has full rank. 
For this reason it is necessary to provide independent points in the initial design. 
This is usually the case, if $d+1$ linearly independent points are provided. 
The matrix inversion can be done efficiently by using singular value decomposition (SVD) or similar algorithms.

The linear polynomial $p(\vec{x})$ in Eq.~\eqref{eq:RBFfunction} serves the purpose to alleviate the fit of simple linear functions $f()$ which otherwise have to be approximated by superimposing many RBFs in a complicated way. Polynomials of higher order may be used as well. We consider here the option of additional direct squares, where $p(\vec{x})$ in Eq.~\eqref{eq:RBFfunction} is replaced by
\begin{equation}
		p_{sq}(\vec{x}) = p(\vec{x}) + e_1 x_1^2 + \ldots + e_d x_d^2
\label{eq:squares}
\end{equation}
with additional coefficients $\vec{e} = (e_1,\ldots,e_d)^T$. The matrix in Eq.~\eqref{eq:RBF} is extended in a straightforward manner from an 
$(n+d+1)\times(n+d+1)$-matrix to an $(n+2d+1) \times (n+2d+1)$-matrix.

%Constrained Optimization by Radial Basis Function Approximation~(COBRA) was proposed by Regis~\cite{regis2014constrained}. The main idea of this method is to use approximations of the objective function and the constraint functions for saving real function evaluations. This can be especially beneficial for industrial optimization problems where function evaluations are expensive. Internally, COBRA uses radial basis function interpolation for modeling the objective and constraints. In each iteration $k$ the solution $\vec{x}^{(k)}$ to be evaluated on the real function $f$ is the result of an optimization performed on the RBF surrogates consisting of surrogates of the objective and the constraint functions.

%%%%%%%%%%%%%%%%%%%%%%%%%%%%%%%%%%%%%%%%%%%%%%%%%%%%%%%%%%%%%%%%%%%%%%%%%%%%%%%%%%%%%%%%%%%%%%%%%%%%
\subsection{Common pitfalls in surrogate-assisted optimization}

RBF models are very fast to train, even in high dimensions. They often provide good approximation accuracy even when only few training points are given. This makes them ideally suited as surrogate models for high-dimensional optimization problems with a large number of constraints.

There are however some pitfalls which should be avoided to achieve good modeling results for any surrogate-assisted black-box optimization.

\begin{figure}[tbp]
	\centering
		\includegraphics[width=0.95\columnwidth]{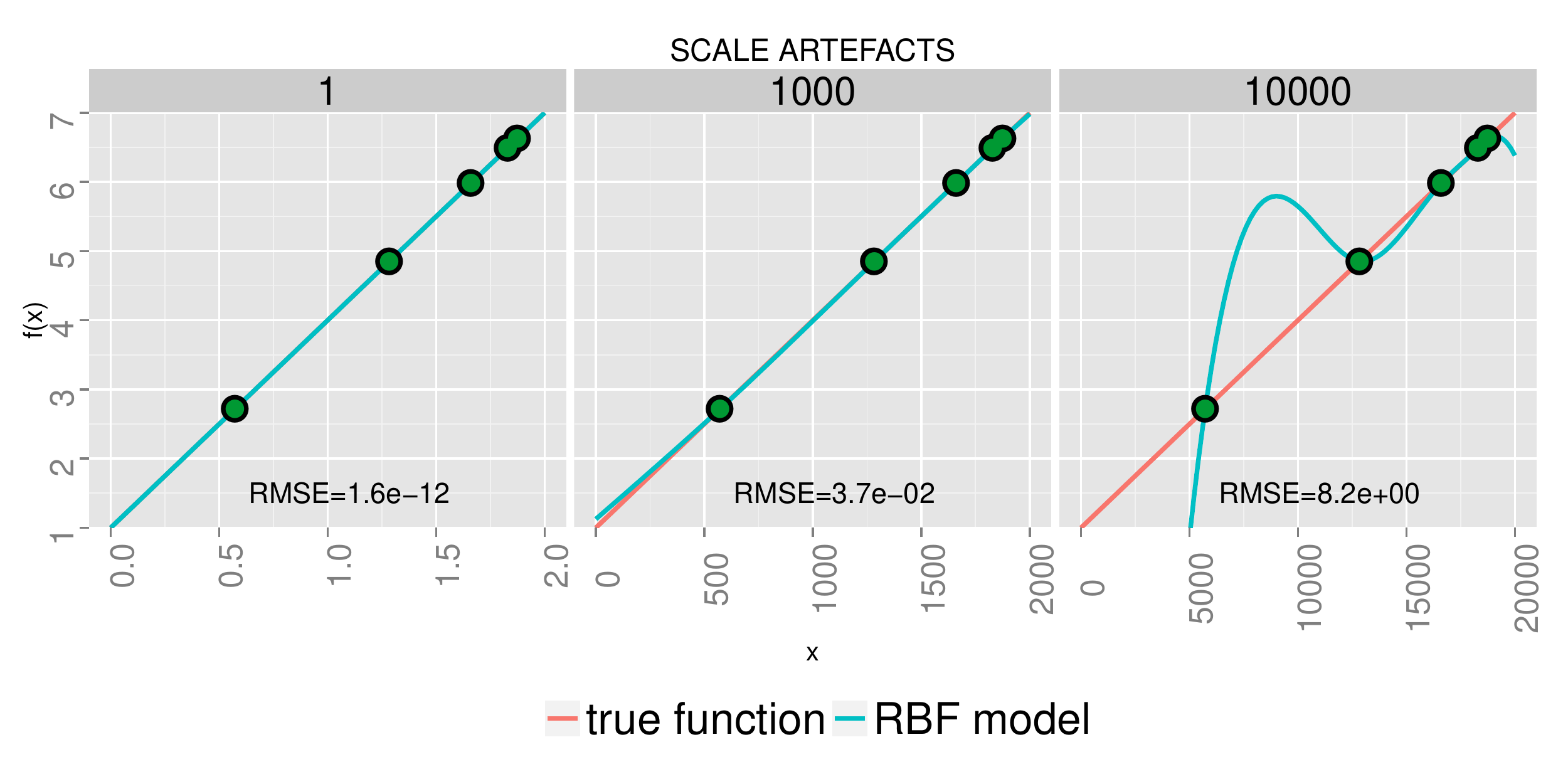}
	\caption{The influence of scaling. From left to right the plots show the RBF model fit for scale $S=1,1000,10000$ (upper facet bar).
	}
	\label{fig:RbfScaling}
\end{figure}

\subsubsection{Rescaling the input space}
\label{Sec:rescale}
If a model is fitted with too large values in input space, a striking failure may occur. Consider the following simple example:
\begin{equation}
		f(x) = 3\frac{x}{S} + 1
\label{eq:RbfLinear}
\end{equation} 
where $x \in [0, 2S]$.  If $S$ is large, the $x$-values (which enter the RBF-model) will be large, although the output produced by Eq.~\eqref{eq:RbfLinear} is exactly the same.  Since the function $f(x)$ to be modeled is exactly linear and the RBF-model contains a linear tail as well, one would expect at first sight a perfect fit (small RMSE) for each surrogate model. But -- as Fig.~\ref{fig:RbfScaling} shows -- this is not the case for large $S$: The fit (based on the same set of five points) is perfect for $S=1$, weaker for $S=1000$, and extremely bad in extrapolation for $S=10000$. 

The reason for this behavior: Large values for $x$ lead to computationally singular (ill-conditioned) coefficient matrices, because the cubic coefficients tend to be many orders of magnitude larger than the coefficients for the linear part. Either the linear equation solver will stop with an error or it produces a result which may have large RMSE, as it is demonstrated in the right plot of Fig.~\ref{fig:RbfScaling}. The solver sets the linear tail of the RBF model to zero in order to avoid numerical instabilities. The RBF model thus attempts to approximate the linear function with a superposition of cubic RBFs. This is bound to fail if the RBF model has to extrapolate beyond the green sample points.

This effect exactly occurs in problem G10, where the objective function is a simple linear function $x_1+x_2+x_3$ and the range for the input dimensions is large and different, e.g. $[100,10000]$ for $x_1$ and $[10,1000]$ for $x_3$. 

The solution to this pitfall is simple: Rescale a given problem in all its input dimensions to a small and identical range, e.g. either to [0,1] or to [-1,1] for all $x_i$.

\begin{figure}[tbp]
	\centering
		\includegraphics[width=0.85\columnwidth]{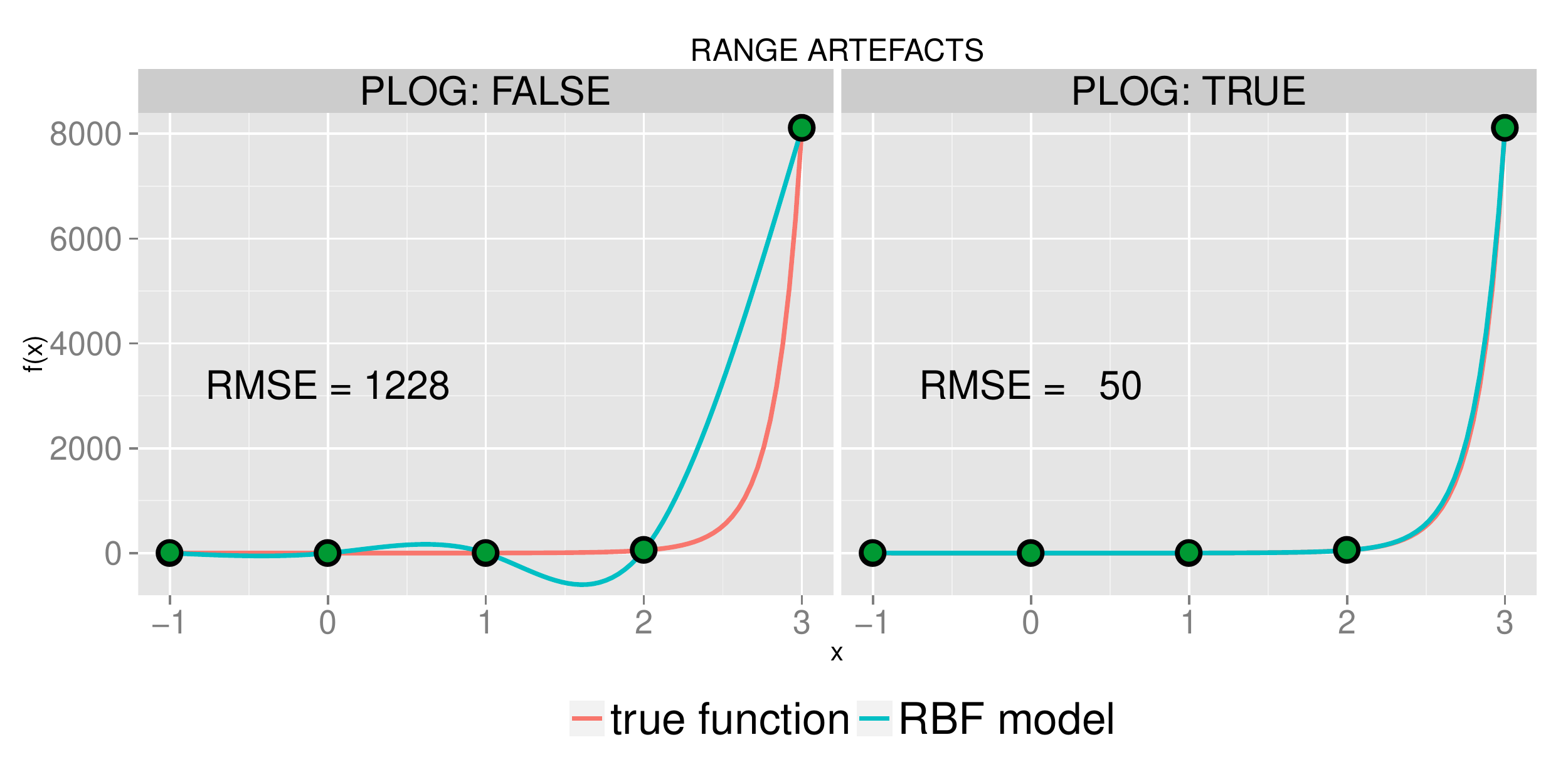}
	\caption{The influence of large output ranges. Left: Fitting the original function with a cubic RBF model. Right: Fitting the $plog$-transformed function with an RBF-model and transforming the fit back to original space with $plog^{-1}$.
	}
	\label{fig:RbfRange}
\end{figure}

\subsubsection{Logarithmic transform for large output ranges}
\label{Sec:plogPitfall}

Another pitfall are large output ranges in objective or constraint functions. As an example consider the function
\begin{equation}
		f(x) = e^{x^2}
\label{eq:RbfExp2}
\end{equation} 
which has small values $<10$ in the interval [-1,1] around its minimum, but quickly grows to large values above 8000 at $x=3$. If we fit the original function with an RBF model using the green sample points shown in Fig.~\ref{fig:RbfRange}, we see in the left plot an oscillating behavior in the RBF function. This results in a large RMSE (approximation error). 

The reason is that the RBF model tries to avoid large slopes. Instead the fitted model is similar to a spline function. Therefore it is a useful remedy to apply a logarithmic transform which puts the output into a smaller range and results in smaller slopes. Regis and Shoemaker ~\cite{regisShoe2012} define the function
\begin{equation}
		\mbox{plog}(y) = \begin{cases}
								 +\ln( 1+y) & \mbox{\quad if \quad} y \geq 0 \\
								 -\ln( 1-y) & \mbox{\quad if \quad} y   <  0
							\end{cases}
\label{eq:plog}
\end{equation} 
which has -- in contrast to the plain logarithm -- no singularities and is strictly monotonous for all $y \in \mathbb{R}$. The RBF model can perfectly fit the $plog$-transformed function. Afterward we transform the fit with $plog^{-1}$ back to the original space and the back-transform takes care of the large slopes. As a result we get a much smaller approximation error RMSE in the original space, as the right-hand side of Fig.~\ref{fig:RbfRange} shows. 

We will apply the $plog$-transform only to \replacedOK{functions with steep slopes in our surrogate-assisted optimization SACOBRA. For functions with flat or constant slope (e.~g. linear functions) our experiments have shown that -- due to the nonlinear nature of $plog$ -- the RBF approximation for  $plog(f)$ is less accurate.}{ large-range functions in our surrogate-assisted optimization SACOBRA.~\SB{This sentence is not any longer true, or am I wrong?} For small-range functions experiments have shown that the $plog$-transform can make it harder for the internal optimizer to locate minima in the surrogate models successfully.}
\subsection{COBRA-R}
\label{Sec:Cobra}

The COBRA algorithm has been developed by Regis~\cite{regis2014constrained} with the aim to solve constrained optimization tasks with severely limited budgets. The main idea of COBRA is to do most of the costly optimization on surrogate models (RBF models, both for the objective function $f$ and the constraint functions $g_i$). We reimplemented this algorithm in \texttt{R}~\cite{Rproject13} with a few small modifications. We give a short review of this algorithm COBRA-R in the following.

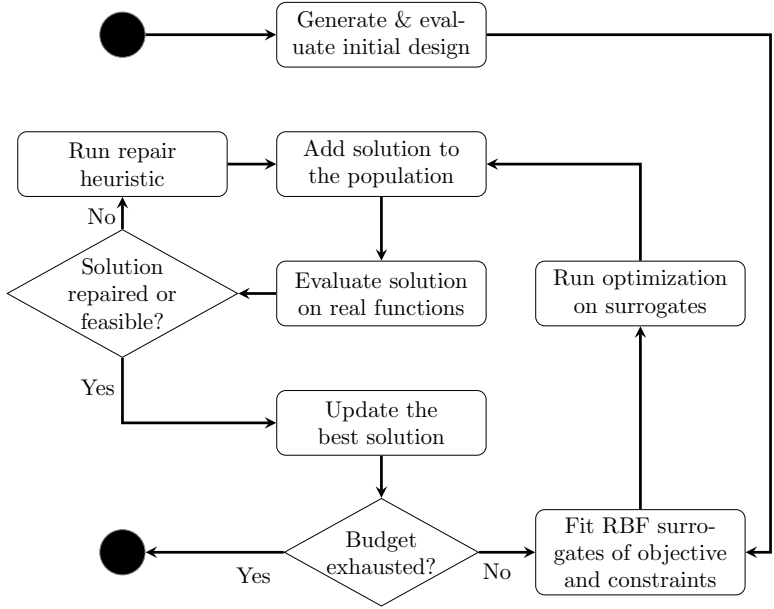
\begin{figure}[tb]
\centerline{
	\resizebox{0.8\textwidth}{!}{
		\input{cobraFlc}
	}
}
\caption{COBRA-R flowchart}%
\label{fig:cobraFlc}%
\end{figure}

COBRA-R starts by generating an initial population $P$ with $n_{init}$ points (i.~e. a random initial design\footnote{usually a latin hypercube sampling (LHS)}, see Fig.~\ref{fig:cobraFlc})  to build the first set of surrogate models . The minimum number of points is $n_{init} =d+1$, but usually a larger choice  $n_{init} =3 d$ gives better results. 

Until the budget is exhausted, the following steps are iterated on the current population $P=\{\vec{x}_1,\ldots,\vec{x}_n\}$: 
%The surrogate models are trained on the current set $P$. 
The constrained optimization problem is executed by optimizing \textit{on the surrogate functions}: 
That is, the true functions $f, g_{1},\ldots,g_{m}$ are approximated with RBF surrogate models $s_{0}^{(n)},s_{1}^{(n)},\ldots,\allowbreak s_{m}^{(n)}$, given the $n$ points in the current population $P$.  In each iteration the COBRA-R algorithm solves with any standard constrained optimizer\footnote{Regis~\protect\cite{regis2014constrained} uses \textsc{MATLAB}'s \textsc{fmincon}, an interior-point optimizer, which is not available in the \texttt{R} environment. In COBRA-R we use mostly Powell's COBYLA, but other constrained optimizer like ISRES are implemented in our \texttt{R}-package \protect\href{https://cran.r-project.org/web/packages/SACOBRA}{\url{https://cran.r-project.org/web/packages/SACOBRA}} as well. } 
the constrained surrogate subproblem  
\begin{eqnarray}
 \text{Minimize}    &  s_{0}^{(n)}(\vec{x}) & \label{eq:sub} \\
 \text{subject to}  &  \vec{x} \in [\vec{a},\vec{b}] \subset \mathbb{R}^d, & \nonumber \\
										&  s_{i}^{(n)}(\vec{x})+\epsilon_{}^{(n)} \leq 0, \quad &i=1,2,\ldots,m  \nonumber\\
					 \text{}  &  \rho_{n}-||\vec{x}-\vec{x}_{j}|| \leq 0 , \quad &j=1,\ldots,n.    \nonumber
\end{eqnarray}
Compared to the original problem in Eq.~\eqref{eq:minprob} this subproblem uses surrogates and it contains two new elements $\epsilon^{(n)}$ and $\rho_{n}$ which are explained in the next subsections. Before going into these details we finish the description of the main loop: The optimizer returns a new solution $\vec{x}_{n+1}$.
If $\vec{x}_{n+1}$ is not feasible, a repair algorithm \textsc{RI2} described in our previous work~\cite{Koch2015a} tries to replace it with a feasible solution in the vicinity.\footnote{\textsc{RI2} is only rarely invoked on the G-problem benchmark but more often in the MOPTA08 case.} In any case, the new solution $\vec{x}_{n+1}$ is evaluated on the true functions $f, g_{1},\ldots,g_{m}$ . It is compared to the best feasible solution found so far and replaces it if better. The new solution $\vec{x}_{n+1}$ is added to the population $P=\{\vec{x}_1,\ldots,\vec{x}_{n+1}\}$ and the next iteration starts with incremented $n$.

\subsubsection{Distance requirement cycle}
COBRA~\cite{regis2014constrained} applies a distance requirement factor which determines how close the next solution $\vec{x}_{n+1} \in \mathbb{R}^d$ is allowed to be to all previous ones. The idea is to avoid frequent updates in the neighborhood of the current best solution. The distance requirement can be passed by the user as external parameter vector $\Xi = \langle \xi^{(1)},\xi^{(2)},\ldots,\xi^{(\kappa)} \rangle$ with $\xi^{(i)} \in \mathbb{R}^{\geq 0}$. In each iteration $n$, COBRA selects the next element $\rho_n = \xi^{(i)}$ of $\Xi$ and 
%measures the distance between the proposed infill solution and the selected element. 
adds the constraints $||\vec{x} - \vec{x}_j|| \geq \rho_n, \quad j=1,...,n$ to the set of constraints. This measures the distance between the proposed infill solution and all $n$ previous infill points. 
The distance requirement cycle (DRC) is a clever idea, since small elements in $\Xi$ lead to more exploitation of the search space, while larger elements lead to more exploration. 
If the last element of $\Xi$ is reached, the selection starts with the first element again. 
% /WK/ I think the sentence commented out below is plain wrong. I tried to reformulate this section. 
%
%If the distance of the previous solution and the proposed infill point is smaller than the element $\xi^{(i)}$, the infill point is discarded. 
%
The size of $\Xi$  and its single components can be arbitrarily chosen.

\subsubsection{Uncertainty of constraint predictions}
COBRA~\cite{regis2014constrained} aims at finding feasible solutions by extensive search on the surrogate functions. However, as the RBF models are probably not exact, especially in the initial phase of the search, a factor $\epsilon^{(n)}$ is used to handle wrong predictions of the constraint surrogates. Starting with $\epsilon^{(n)}=0.005\cdot l$, where $l$ is the diameter of the search space,  
a point $\vec{x}$ is said to be \textit{feasible in iteration $n$} if 
\begin{equation}
s_i^{(n)}(\vec{x}) + \epsilon^{(n)} \leq 0 \quad \forall \quad i=1,\ldots,m
\end{equation}
holds.
That is, we tighten the constraints by adding the factor $\epsilon^{(n)}$ which is adapted during the search. The $\epsilon^{(n)}$-adaptation is done by counting the feasible and infeasible infill points $C_{feas}$ and $C_{infeas}$ over the last iterations. When the number of these counters reaches the threshold for feasible or infeasible solutions, $T_{feas}$ or $T_{infeas}$, respectively, we divide or double $\epsilon^{(n)}$ by $2$ (up to a given maximum $\epsilon_{max}$). When $\epsilon^{(n)}$ is decreased, solutions are allowed to move closer to the real constraint boundaries (the imaginary boundary is relaxed), since the last $T_{feas}$ infill points were feasible. Otherwise, when no feasible infill point is found for $T_{infeas}$ iterations, $\epsilon^{(n)}$ is increased in order to keep the points further away from the real constraint boundary.

\subsubsection{Differences COBRA vs. COBRA-R} %~\\
%\WK{If we need too much words to describe this correctly, we may as well consider to move this to an appendix. Let us leave it for the moment here and decide later.}
Although COBRA and COBRA-R are sharing many common principles, there are several differences which can lead to different results on identical problems. The main differences between COBRA~\cite{regis2014constrained} and COBRA-R~\cite{Koch14a} are listed as follows:
\begin{itemize}
	%\item COBRA-R is implemented in \texttt{R}, but COBRA is implemented in \texttt{MATLAB}. 
	%\item COBRA is flexible in terms of selecting the internal optimizer. Any constrained or unconstrained optimization technique can be easily embedded in the COBRA-R framework and be utilized. We used, the available optimization packages in \texttt{R} like COBYLA, ISRES, HJKB, NMKB.
	\item 
	%\emph{Implementation environment}: 
	COBRA is implemented in \textsc{MATLAB} while COBRA-R is implemented in \texttt{R}. 
  \item \emph{Internal optimizer:} COBRA uses \textsc{MATLAB}'s \textsc{fmincon}, an interior-point optimizer, COBRA-R uses COBYLA.\footnote{Other optimizers like ISRES and unconstrained optimizers with penalty are also available in COBRA-R / SACOBRA package, but not used in this paper.}
	%\item \emph{Initialization} 
	\item \emph{Skipping phase 1:} COBRA has an additional phase 1 for searching the first feasible point.\footnote{ 
	Phase 1 uses an objective function which rewards constraint fulfillment. We implemented this in COBRA-R as well but found it to be unnecessary for our problems. In this paper, COBRA-R always skips phase 1 and directly proceeds with phase 2 even if no feasible solution is found in the initialization phase.}
	\item \emph{Repair algorithm:} COBRA-R has an additional repair algorithm  RI2~\cite{Koch2015a}.  
	\item \emph{Rescaling the input space}: COBRA rescales each dimension to [0,1], COBRA-R in its initial form~\cite{Koch14a} does not. See Sec.~\ref{Sec:rescaleSA} for further remarks on rescaling.
\end{itemize}

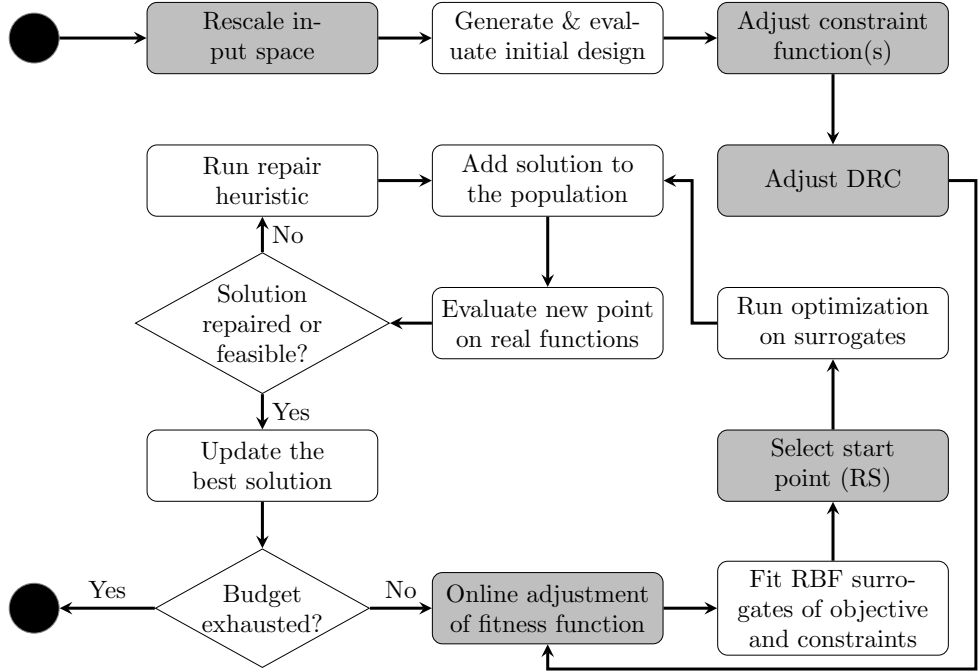
\begin{figure}[tb]
\centerline{
	\resizebox{1.0\textwidth}{!}{
		\input{sacobra2Flc}
	}
}
\caption{SACOBRA flowchart}%
\label{fig:sacobra2Flc}%
\end{figure}

%%%%%%%%%%%%%%%%%%%%%%%%%%%%%%%%%%%%%%%%%%%%%%%%%%%%%%%%%%%%%%%%%%%%%%%%%%%%%%%%%%%%%%%%%%%%%%%%%%%%
\subsection{SACOBRA}
\label{Sec:Sacobra}

%\WK{I think, the presentation of both SACOBRA1 and SACOBRA2 as flowchart and as algorithm is too much for the reader. We should probably concentrate on SACOBRA2 as the main SACOBRA algorithm. I shift the SACOBRA1 algo to the end of the document for reference, but we probably leave it out later. 
%}

\input{ASOC-sacobra2Alg}

COBRA and COBRA-R achieve good results on most of the G-problems and on MOPTA08 as studies from Regis~\cite{regis2014constrained} and our previous work~\cite{Koch2015a} have shown. However, it was necessary in both papers to carefully adjust the parameters of the algorithm to each problem and sometimes even to modify the problem by applying a $plog$-transform (Eq.~\eqref{eq:plog}) to the objective function \addOK{or linear transformations to the constraints} or by rescaling the input space. In real black-box optimization all 
%these adjustments are not possible or make the usage of the optimization scheme too complicated.
these adjustments would probably require knowledge of the problem or several executions of the optimization code otherwise.

It is the main contribution of the current paper to present with SACOBRA (Self-Adjusting COBRA) an enhanced COBRA algorithm which has no needs for manual adjustment to the problem at hand. Instead, SACOBRA extracts during its execution information about the specific problem (either after the initialization phase or online during iterations) and takes internally appropriate measures to adjust its parameters or to transform functions. 

We present in Fig.~\ref{fig:sacobra2Flc} the flowchart of SACOBRA where the five new elements compared to COBRA-R are highlighted as gray boxes. The complete SACOBRA algorithm is presented in detail in Algorithm~\ref{alg:sacobra2} -- \ref{alg:rStart2}. We describe in the following the five new elements in the order of their appearance: 

\subsubsection{Rescaling the input space}
\label{Sec:rescaleSA}
 The input vector $\vec{x}$ is element-wise rescaled to $[-1,+1]$. This is done before the initialization phase. It helps to have a better exploration all over the search space because all dimensions are treated the same. More importantly, it avoids numerical instabilities caused by high values of $\vec{x}$ as shown in Sec.~\ref{Sec:rescale}.
	%modeling issues in the case of having highly varied ranges in different dimensions.  

\subsubsection{Adjusting constraint function(s) (aCF)} 
aCF is done by normalizing the range of constraint functions for each problem. The range 
$\widehat{GR}_i$ for the $i$th constraint is estimated from the initial population in Algorithm~\ref{alg:sacobra2Adjust}. Normalizing each $\widehat{GR}_i$ by the average constraint range 
\begin{equation}
		\mbox{avg}\left( \widehat{GR}_i \right) = \frac{1}{m}\sum_i{\widehat{GR}_i}
\label{eq:widehat_GR}
\end{equation}
helps to shift the range of all constraints as little as possible.  Including this step aCF boosts up the optimization performance because all constraints operate now in a similar range. 

\subsubsection{Adjusting DRC parameter (aDRC)} 
aDRC is done after the initialization phase. 
	Our experimental analysis showed that large DRC values can be harmful for problems with a very steep objective function, because a large move in the input space yields a very large change in the output space. This may spoil the RBF model in a sense similar to Sec.~\ref{Sec:plogPitfall} and lead in consequence to large approximation errors. 
	Therefore, we developed an automatic DRC adjustment which selects the appropriate DRC set according to the information extracted after the initialization phase. \addOK{Function \textsc{AnalyzePlogEffect} in Algorithm~\ref{alg:sacobra2Adjust} selects the 'small' DRC $\Xi_s$ if the estimated objective function range $\widehat{FR}$ is larger than a threshold, otherwise it selects the 'large' DRC $\Xi_l$.}
	%\WK{It would be nice if we could give a more reasonable argument  for aDRC than just 'experimental analysis showed ...'.}\SB{see Fig.~\ref{fig:DRC1} and \ref{fig:DRC2}}
	%\SB{I tried to change this a little bit} 
	%\WK{Perhaps we can show that with a large DRC element the approximation error (say, at the true optimum, to have a common reference point) increases significantly.}
	
\subsubsection{Random start algorithm (\textsc{RS})} 
	Normally COBRA starts optimization from the \replacedOK{current best point}{best point found so far}. With RS (Algorithm~\ref{alg:rStart2}), the optimization starts from a random point in the search space with a certain constant probability $p_{1}$.
%or after $N_{max}/10$ unsuccessful iterations in a row. 
\addOK{If the rate of feasibile individuals in the population $P$ drops below 5\%  then we replace $p_1$ with a larger probability $p_{2}$.} \textsc{RS} is especially beneficial when the search gets stuck in local optima or when it gets stuck in a region where no feasible point can be found.

\subsubsection{Online adjustment of fitness function (aFF)} 
\label{Sec:aFF}
%\SB{I changed this part a little bit}	
Our analysis in Sec.~\ref{Sec:plogPitfall} has shown that a fitness function $f$ with steep slopes poses a problem for RBF approximation. For some problems, modeling $plog(f)$ instead of $f$ and transforming the RBF result back with $plog^{-1}$ boosts up the optimization performance significantly.
%It is better to  model $plog(f)$ and transform the RBF result back with $plog^{-1}$. 
On the other hand, our tests have shown that the $plog$-transform is harmful for some other problems. 
Therefore, a careful \replacedOK{decision whether to use $plog$ or not should be made}{ selection should be done}.

 The idea of \replacedOK{our online adjustment algorithm}{ aFF} (Algorithm~\ref{alg:sacobra2Adjust}, function \textsc{AnalyzePlogEffect}) is the following:
Given the population $P$, we build RBFs for $f$ and $plog(f)$, take a new point $\vec{x}_{new}$ not yet added to $P$, and calculate the ratio of approximation errors on $\vec{x}_{new}$ (line 15 of Algorithm~\ref{alg:sacobra2Adjust}). We do this in every $k$th iteration (usually $k=10$) and collect these ratios in a set $E$. 
If %the median $Q$ of  set $E$ 
\begin{equation}
			Q = \log_{10} \left( \med(E) \right)
\label{eq:Qvalue}
\end{equation}
is above 0, then the RBF for $plog(f)$ is better in the majority of the cases. Otherwise, the RBF on $f$ is better.\footnote{Our experimental analysis on the G-problem test suite will show (Sec.~\protect\ref{Sec:fitfuncadjust}) that a threshold 1 is slightly more robust than 0. We use this threshold 1 in step 11 of Algorithm~\ref{alg:sacobra2}, but the difference to threshold 0 is only marginal.}

%Function \textsc{AdjustFitnessFunction}  is now commented out, it is directly coded in Step 11 of \ref{alg:sacobra2}.
Step 11 of Algorithm~\ref{alg:sacobra2} decides on the basis of this criterion $Q$ which function $\widetilde{f}$ is used as RBF surrogate in the optimization step. Note that the decision for $\widetilde{f}$ taken in earlier iterations can be revoked in later iterations, if the majority of the elements in $E$ shows
%\SB{show or shows?} \WK{shows, it refers to noun majority}
that now the other choice is more promising. 

\vspace{0.7cm}
This completes the description of our SACOBRA algorithm. 
SACOBRA is available as open-source \texttt{R}-package from CRAN.\footnote{\protect
\href{https://cran.r-project.org/web/packages/SACOBRA}{\url{https://cran.r-project.org/web/packages/SACOBRA}}
}
 
%%%%%%%%%%%%%%%%%%%%%%%%%%%%%%%%%%%%%%%%%%%%%%%%%%%%%%%%%%%%%%%%%%%%%%%%%%%%%%%%%%%%%%%%%%%%%%%%%%%%
%\subsection{Performance and Data Profiles}
\subsection{Performance Measures}
\label{Sec:PerfMeasure}
\par{In many papers on optimization the strength of an optimization technique is measured by comparing the final solution achieved by different algorithms~\cite{runarsson2000stochastic}. This approach only provides the information about the quality of the results and neglects  the speed of convergence which is a very important measure for expensive optimization problems. Comparing the convergence curve over time (number of function evaluations) is also one of the common benchmarking approaches~\cite{regis2014constrained}. Although a convergence curve provides good information about the speed of convergence and the final quality of the optimization result, it can be used to compare performance of several algorithms only on \textit{one} problem. It is often interesting to compare the overall capability of a technique on solving a group of problems. The data and performance profiles developed by Mor\'e and Wild~\cite{WildMore2009} are a good approach to analyze \addOK{the} performance of any optimization algorithm on a whole test suite and are now used frequently in the optimization literature~\cite{brockhoff2015,regis2015trust}}. 
%\WK{Why should we cite 'no free lunch' here?}\SB{It was wrong, omitted}

\subsubsection{Performance Profiles} 
Performance profiles are defined with the help of the performance ratio
\begin{equation}
r_{p,s}=\frac{t_{p,s}}{\min\limits_{\forall s' \in\,  \mathbb{S}}\{t_{p,s'}\}}, \qquad p \in \mathbb{P}
\end{equation}
where $\mathbb{P}$ is a set of problems, $\mathbb{S}$ is a set of solvers and $t_{p,s}$ is the number of iterations solver $s \in \mathbb{S}$ requires to solve problem $p \in \mathbb{P}$. 
A problem is said to be \textit{solved} when a feasible objective value $f(x)$ is found which is not more than $\tau$ larger than the best objective $f_L$ determined by any solver in $S$:
\begin{equation}
f(x)-f_L \leq \tau
\label{eq:convTest}
\end{equation}
We use $\tau=0.05$ for all our experiments below. 
%\WK{I added this sentence. Is this correct? Is it what is named "tolerance" in the title of the data profile plots?}\SB{correct.}
Smaller values are more desirable for \addOK{the} performance ratio $r_{p,s}$. When using the best solver $s$ to solve problem $p$ then $r_{p,s}=1$. If a solver $s$ cannot solve problem $p$ the performance ratio is set to infinity. The performance profile $\rho_{s}$ is now defined as a function of the steerable performance factor $\alpha$:
\begin{equation}
\rho_{s}(\alpha)=\frac{1}{|\mathbb{P}|}|\{p \in \mathbb{P} : r_{p,s} \leq \alpha\}|.
\end{equation}

\subsubsection{Data Profiles}
Data profiles are appropriate for evaluating optimization algorithms on expensive problems. They are defined as
\begin{equation}
		d_{s}(\alpha)=\frac{1}{| \mathbb{P}|}|\{p \in  \mathbb{P} : \frac{t_{p,s}}{d_{p}+1} \leq \alpha\}|,
\end{equation}
%where $\mathbb{P}$ is a set of problems, $\mathbb{S}$ is a set solvers and $t_{p,s}$ is the number of iterations solver $s \in \mathbb{S}$ required to solve problem $p \in \mathbb{P}$.
with $\mathbb{P}, \mathbb{S}$ and $t_{p,s}$ defined as above and $d_p$ as the dimension of problem $p$.

We prefer data profiles over performance profiles, because the performance factor $\alpha$ has a more intuitive meaning for data profiles: If we allow for each problem with dimension $d_p$ a budget of $B_{\alpha}=\alpha(d_p+1)$ function evaluations, then the value $d_{s}(\alpha)$ can be interpreted as the fraction of problems which solver $s$ can solve within this budget $B_{\alpha}$.

%\WK{We need in Performance Profile an equation how to come from $r_{p,s}$ to the performance profile curve. }\SB{done}
%\WK{And for both curves a definition how 'a problem is solved' is defined.}\SB{done}

%%%%%%%%%%%%%%%%%%%%%%%%%%%%%%%%%%%%%%%%%%%%%%%%%%%%%%%%%%%%%%%%%%%%%%%%%%%%%%%%%%%%%%%%%%%%%%%%%%%%
\section{Experiments}
%%%%%%%%%%%%%%%%%%%%%%%%%%%%%%%%%%%%%%%%%%%%%%%%%%%%%%%%%%%%%%%%%%%%%%%%%%%%%%%%%%%%%%%%%%%%%%%%%%%%
\label{Sec:Experiments}

%%%%%%%%%%%%%%%%%%%%%%%%%%%%%%%%%%%%%%%%%%%%%%%%%%%%%%%%%%%%%%%%%%%%%%%%%%%%%%%%%%%%%%%%%%%%%%%%%%%%
\subsection{Experimental Setup}
\label{Sec:ExpSetup}

\begin{table}[tbp]
\centering
\caption[Characteristics of the G-functions]{Characteristics of the G-functions: $d$: dimension, type of fitness function, $\rho^{*}$: feasibility rate (\%) after changing equality constraints to inequality constraints, $FR$: range of the fitness values, $GR$: ratio of largest to smallest constraint range, LI: number of linear inequalities, NI: number of nonlinear inequalities, NE: number of nonlinear equalities, $a$: number of constraints active at the optimum.}
\label{tab:GprobFeatures}
%	\resizebox{0.8\textwidth}{!}{%
\begin{tabular}{lrrrrrrrrr}
\toprule
Fct. 	& $d$		& type			  & $\rho^{*}$		 & $FR$					  &$GR$	& ~LI 	& NI	& NE    &$\,~a$\\
\midrule
G01 	& 13	& quadratic 	  & 0.0003\%			 & 298.14					&1.969  & 9		&	0		&	0		& 6 \\
G02 	& 10	& nonlinear 	 	& 99.997\%			 &	0.57	  			&2.632  &	1		&	1		&	0	  & 1 \\
G03 	& 20	& nonlinear 	  & 0.0000\%			 &	92684985979.23&1.000  &	0		&	0		&	1	  & 1 \\
G04		& 5		& quadratic 		& 26.9217\%			 &	9832.45				&2.161	&	0		&	6		&	0		& 2 \\
G05		& 4		& nonlinear 	  & 0.0919\%			 &	8863.69				&1788.74&	2		&	0		&	3	  & 3 \\
G06		& 2		& nonlinear 	  & 0.0072\%			 &	1246828.23		&1.010	&	0		&	2		&	0	  & 2 \\
G07		& 10	& quadratic 	  & 0.0000\%			 &	5928.19				&12.671	&	3		&	5		&	0	  & 6 \\
G08		& 2		& nonlinear 	  & 0.8751\%			 &	1821.61				&2.393	&	0		&	2		&	0	  & 0 \\
G09		& 7		& nonlinear 	  & 0.5207\%			 &	10013016.18		&25.05	&	0		&	4		&	0	  & 2 \\
G10		& 8		& linear    	  & 0.0008\%			 &	27610.89			&3842702&	3		&	3		&	0	  & 3 \\
G11		& 2		& linear    	  & 66.7240\%			 &	4.99			    &1.000	&	0		&	0		&	1	  & 1 \\
\bottomrule
\end{tabular}
%}  % resize
\end{table}

We evaluate SACOBRA by using a well-studied test suite of G-problems described in~\cite{Floudas,michalewicz1996evolutionary}. The diversity of the G-problem characteristics makes them a very challenging benchmark for optimization techniques. In Table~\ref{tab:GprobFeatures} we show features of these problems. The features $\rho^{*}$, $FR$ and $GR$ (defined in Table~\ref{tab:GprobFeatures}) are measured by Monte Carlo sampling with $10^6$ points in the search space of each G-problem.

Equality constraints are treated by replacing each equality operator with an inequality operator of the appropriate direction. This approach (same as in Regis' work~\cite{regis2014constrained}) takes as \glqq appropriate direction\grqq\ this side of the equality hyperplane where the objective function increases. 

The MOPTA08 benchmark by Jones~\cite{jones2008large} is a substitute for a high-dimen\-sio\-nal real-world problem encountered in the automotive industry: It is a problem with $d=124$ dimensions and with 68 constraints. 
The problem should be solved within $1860 = 15\cdot d$ function evaluations. This corresponds to one month of computation time on a high-performance computer for the real automotive problem since the real problem requires time-consuming crash-test simulations.

%\paragraph{Initialization}~\\
The COBRA-R optimization framework allows the user to choose between several initialization approaches: Latin hypercube sampling (\textsc{LHS}), \textsc{Biased} and \textsc{Optimized}~\cite{Koch14a}. While \textsc{LHS} initialization is always possible (and is in fact used for all runs of the G-problem benchmark with $n_{init}=3d$), the other algorithms are only possible if a feasible starting point is provided. In COBRA~\cite{regis2014constrained} the initialization is always done randomly by means of Latin hypercube sampling for functions without feasible starting point.

In the case of MOPTA08 a feasible point is known. We use the \textsc{Optimized} initialization approach, where an initial optimization run is started from this feasible point with the Hooke \& Jeeves pattern search algorithm~\cite{hooke1961direct}. This initial run provides a set of $n_{init}=500$ points in the vicinity of the feasible point. This set serves as initial design for MOPTA08.

%\WK{A question comes to my mind: If COBRA can solve MOPTA08 with LHS initial sampling: Why can we not do this with COBRA-R?}
%\SB{I just corrected this sentence in COBRA-R we also use optimized initialization for MOPTA08}

\input{settingTab}
Table~\ref{tab:setting} shows the parameter settings used for COBRA-R and SACOBRA. All G-problems were optimized with exactly the same initial parameter settings. In contrast to that, the COBRA results in Regis~\cite{regis2014constrained} and our previous work~\cite{Koch14a} were obtained by manually activating $plog$ for some G-problems and by manually adjusting constraint factors and other parameters.

%This figure is generated by conCurve.R
% We use SACOBRA50 to generate these figures
% See doc\Notes.d\notes-SB-MONREP\SACOBRA-test.xlsx for in-detail settings of SACOBRA50. The most important settings are:
%   constant restart probability p1=0.175, which is only augmented to p2=0.4 if the feasibility rate in the population drops below 5%;
%   cubic RBF, repairMargin=0.1, onlinePLOG=T.
% The current runs return feasible solutions for all G-functions
%BOTHBANDS=T
%---------------------------------------------------------------------------------------------------------------
\begin{figure}%
\includegraphics[width=\columnwidth]{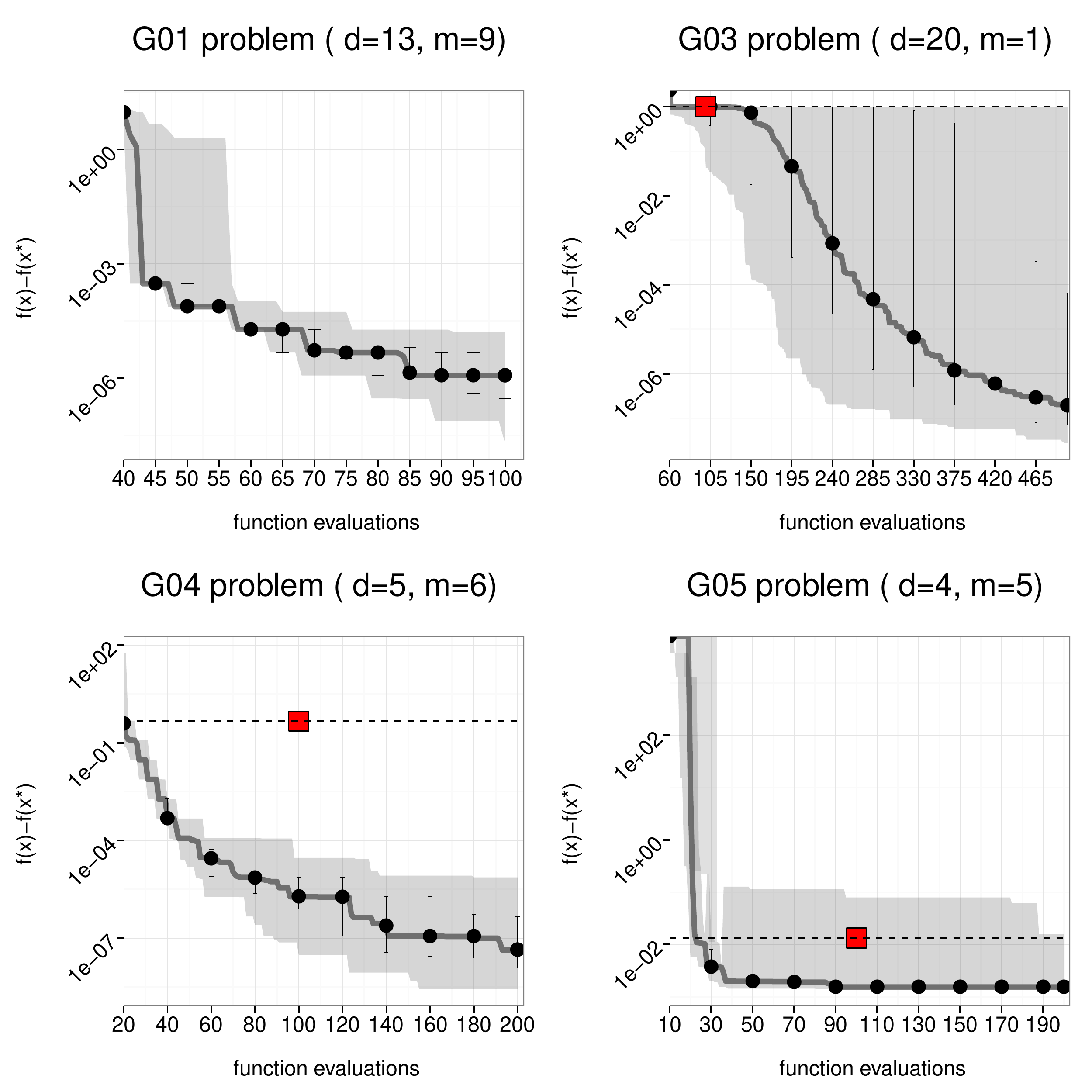}%
\caption{SACOBRA optimization process for G01 -- G05. The gray curve shows the median of the error for 30 independent trials. The error is calculated with respect to the true minimum $f(x^*)$. The gray shade around the median is showing the worst and the best error. The error bars mark the 25\% and 75\% quartile. The red square is the result reported by Regis~\cite{regis2014constrained} after 100 iterations. 
}%
\label{fig:conv1}%
\end{figure}
%---------------------------------------------------------------------------------------------------------------

%This figure is generated by conCurve.R
% We use SACOBRA50 to generate these figures
% See doc\Notes.d\notes-SB-MONREP\SACOBRA-test.xlsx for in-detail settings of SACOBRA50. 
% The current runs return feasible solutions for all G-functions
%BOTHBANDS=T
%---------------------------------------------------------------------------------------------------------------
\begin{figure}%
\includegraphics[width=\columnwidth]{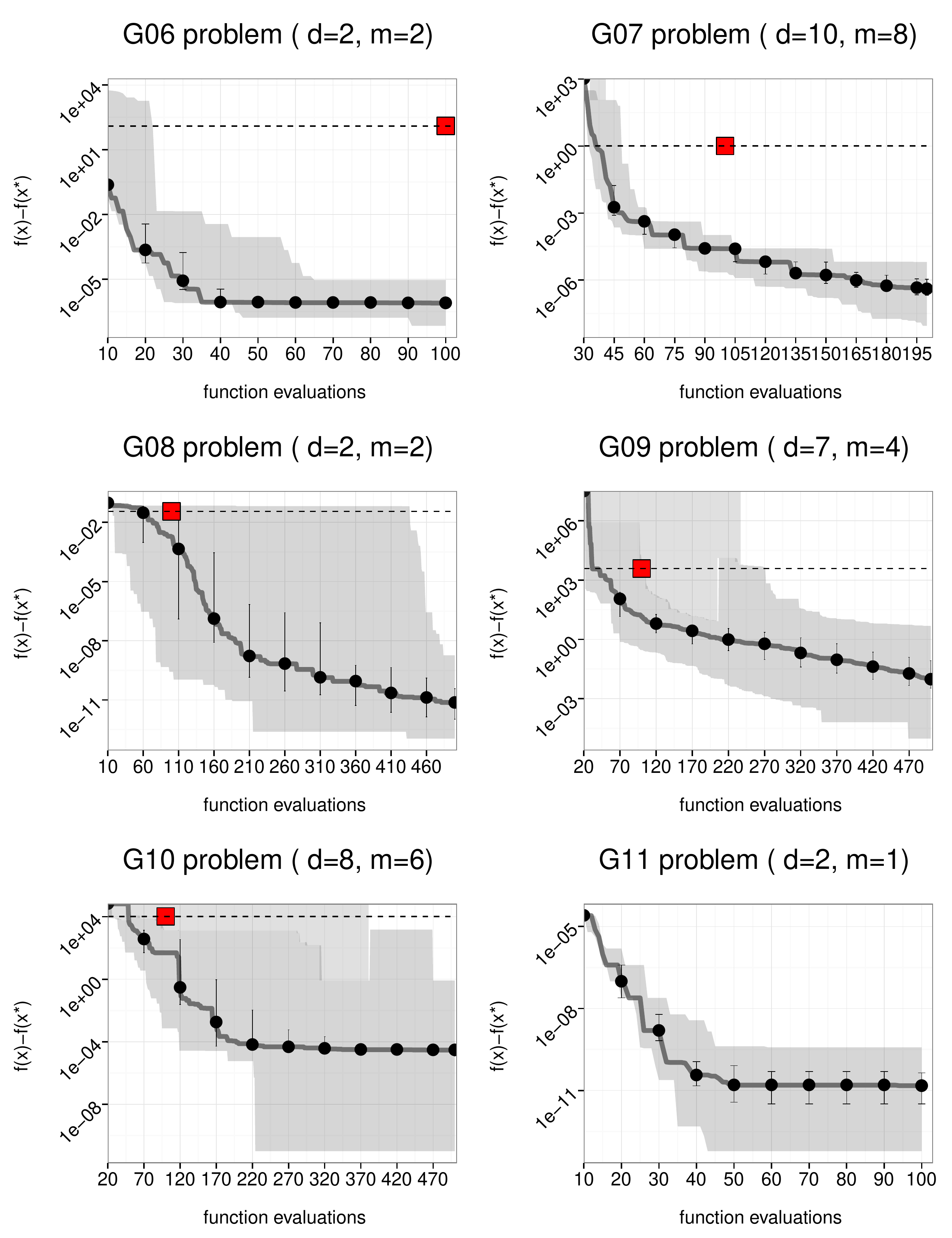}%
\caption{Same as Fig.~\ref{fig:conv1} for G06 -- G11.
}%
\label{fig:conv2}%
\end{figure}
%---------------------------------------------------------------------------------------------------------------

%This figure is generated by conCurve.R
% We use SACOBRA50 to generate these figures
% See doc\Notes.d\notes-SB-MONREP\SACOBRA-test.xlsx for in-detail settings of SACOBRA50. 
% The current runs return feasible solutions for all G-functions
\begin{figure}%
\begin{minipage}[b]{0.5\textwidth}
\includegraphics[width=0.9\columnwidth]{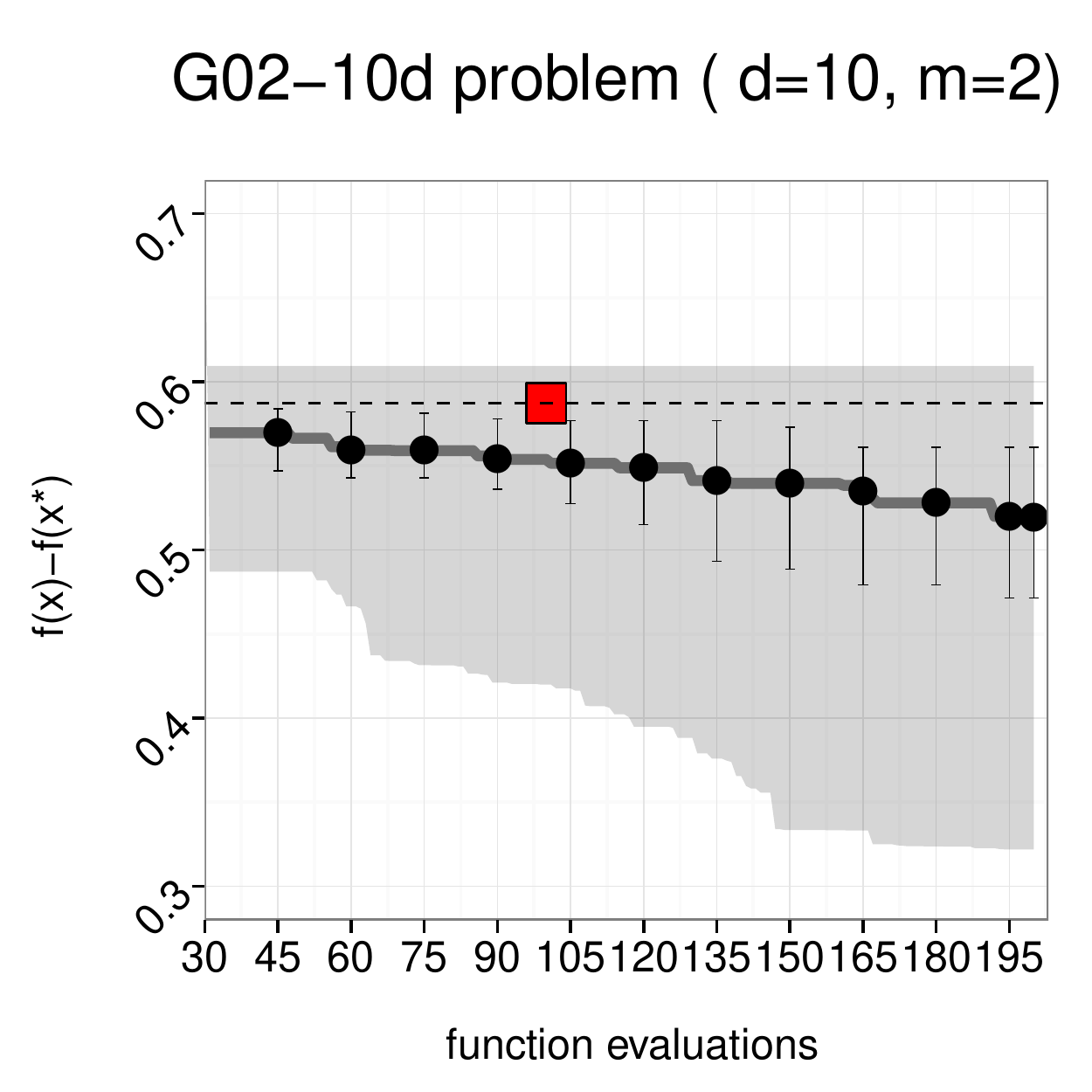}%
\end{minipage}
\hfill
	\begin{minipage}[b]{0.5\textwidth}
\includegraphics[width=0.9\columnwidth]{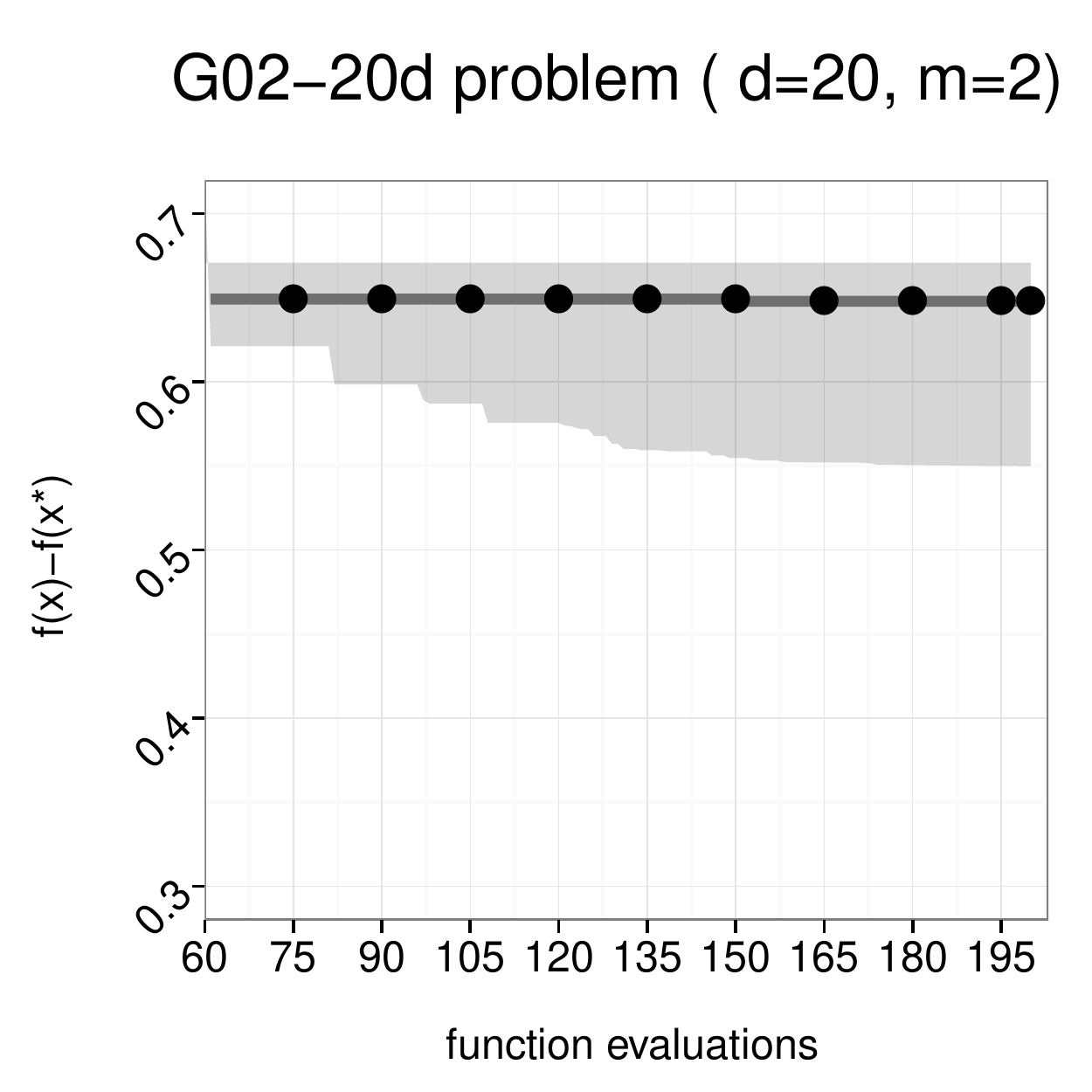}%
	\end{minipage}
%\end{tabular}
\caption{Same as Fig.~\ref{fig:conv1} for G02 in 10 and 20 dimensions.}%
\label{fig:conv-g02}%
\end{figure}

%%%%%%%%%%%%%%%%%%%%%%%%%%%%%%%%%%%%%%%%%%%%%%%%%%%%%%%%%%%%%%%%%%%%%%%%%%%%%%%%%%%%%%%%%%%%%%%%%%%%
\subsection{Convergence Curves}
\label{Sec:convcurves}

Figures~\ref{fig:conv1} -- \ref{fig:conv-g02} show the SACOBRA convergence plots for all G-problems. It is clearly visible that all problems except G02 are solved in the majority of runs, if we define \textit{solved} as a target error below $\tau=0.05$ in comparison to the true optimum.  In some cases (G03, G05, G09, G10) the worst error does not meet the target, but in the other cases it does. In most cases, as indicated by the red squares, there is a clear improvement to Regis' COBRA results~\cite{regis2014constrained}.

%%%%%%%%%%%%%%%%%%%%%%%%%%%%%%%%%%%%%%%%%%%%%%%%%%%%%%%%%%%%%%%%%%%%%%%%%%%%%%%%%%%%%%%%%%%%%%%%%%%%
\subsection{Performance Profiles}
\label{Sec:perfprofile}

%----------------------a plot comparing PF and DF of differnet versions of SACOBRA in absence of different-----------------------------------
%---------------------------sacobra elements like: rescale, RS, aFF, aCF and ...-------------------------------------------------------------
%-----SACOBRA-compariosn figure in color
%generated by SACOBRAcomparison.R with the folowing parameters
%RegisRule=TRUE
%tol<-0.05
%type<-"DP"
%COLORED=TRUE
%SHAPED=TRUE
%methods<-c("SACOBRA50","SACOBRA51","SACOBRA52","SACOBRA53","SACOBRA54","SACOBRA55","SACOBRA32","SACOBRA9")
\begin{figure}[tb]%
\includegraphics[width=0.9\columnwidth]{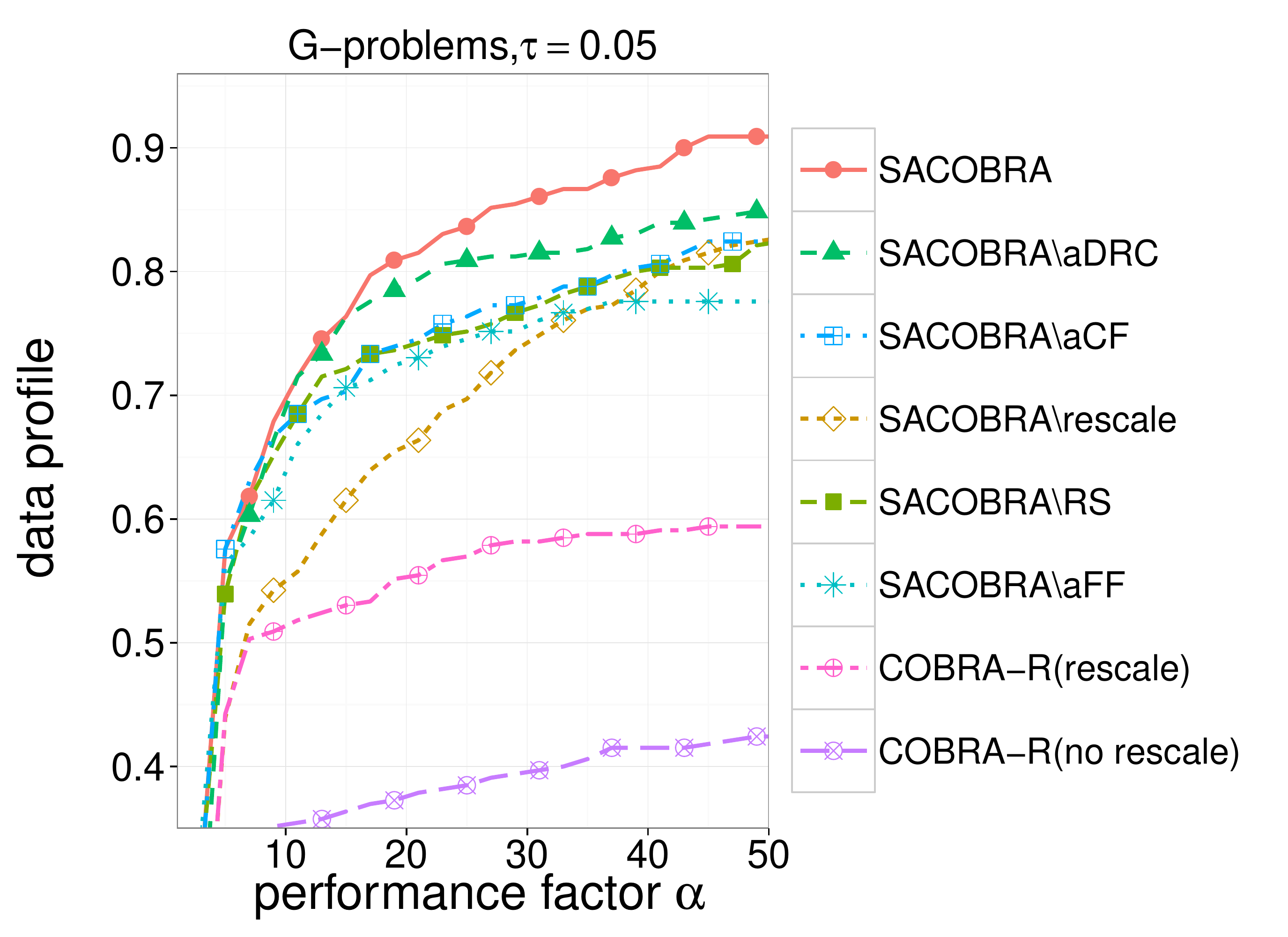}%
\caption{Analyzing the impact of different elements of SACOBRA on the G-problems. Data profile of SACOBRA, SACOBRA$\backslash$rescale (SACOBRA without rescaling the input space), and other \glqq $\backslash$\grqq-algorithms are  with a similar meaning. COBRA-R is the old version of COBRA~\cite{Koch14a}, i.~e. SACOBRA with all adjustment extensions switched off. These algorithms are performed on 330 different problems (11 test problems from G-function suite which are initialized with 30 different initial design points).
%\WK{replace \textit{tolerance} by $\tau$ in the plot title}
}%
\label{fig:SACOBRAcompfinal}%
\end{figure}

Our main result is shown in Fig.\ref{fig:SACOBRAcompfinal}. It shows the data profiles for different SACOBRA variants in comparison with the data profile for COBRA-R. COBRA-R was run with a fixed parameter set.\footnote{In our previous work~\protect\cite{Koch14a,Koch2015a} we reported good results with COBRA-R, but this was with varying parameters and with tedious parameter tuning on each specific G-problem.} We note in passing that other fixed parameter settings for COBRA-R were tested, they were perhaps better on some of the runs but inevitably worse on other runs, so that in the end a similar or slightly worse data profile for COBRA-R would emerge. SACOBRA increases significantly the success rate on the G-problem benchmark suite. 

In addition, we analyze in Fig.~\ref{fig:SACOBRAcompfinal} the effect of the five elements of SACOBRA: The $\backslash$-data profiles present the SACOBRA results when one specific of the five SACOBRA elements is switched off. We see that the strongest effects occur when rescale is switched off (early iterations) or when aFF is switched off (later iterations).  

\begin{figure}[tb]%
\centerline{
\includegraphics[width=0.7\columnwidth,height=0.6\columnwidth]{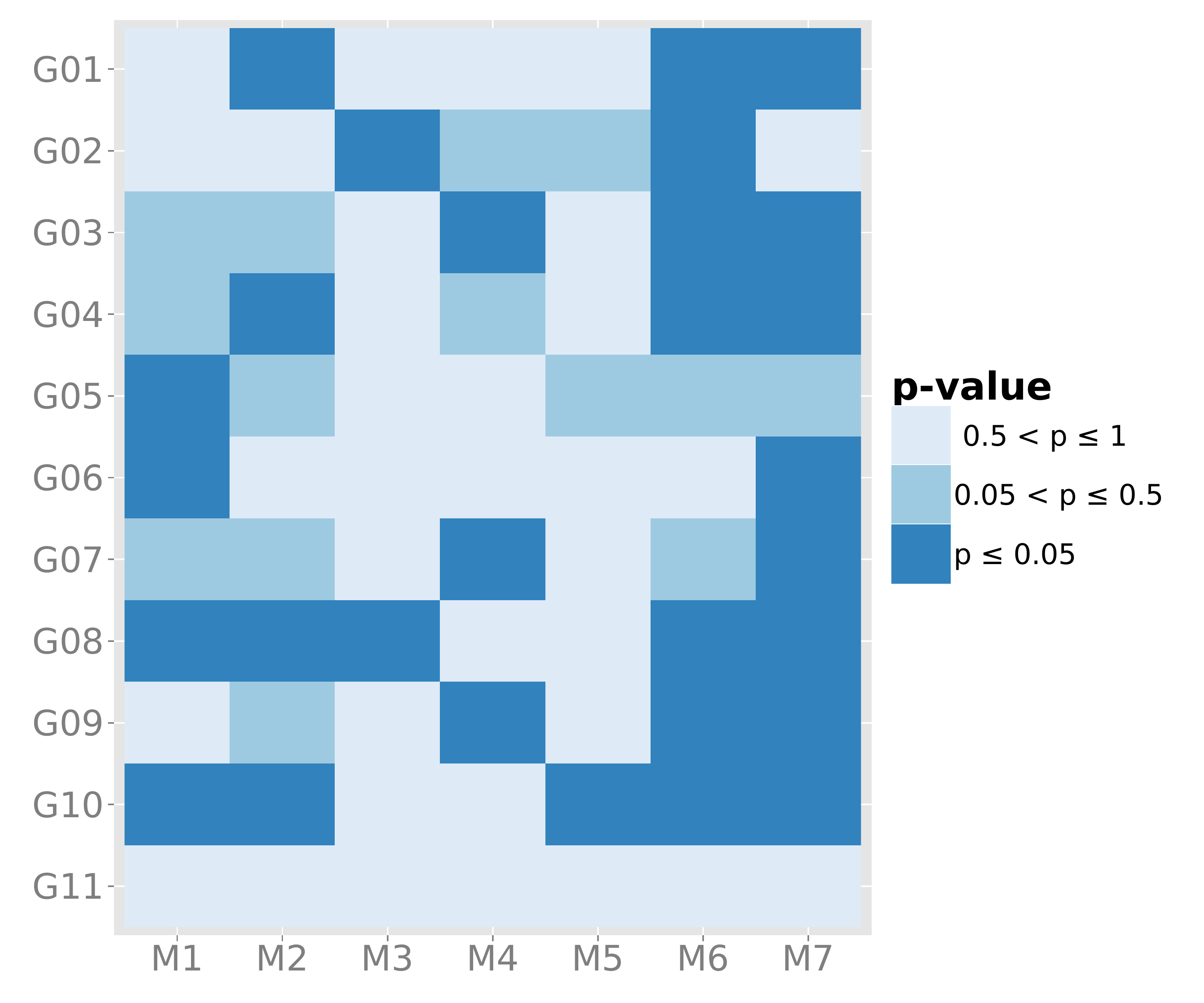}%
%%%%
%%%% generated with Figures.d/wilcox.R
%%%% SS=TRUE, BLUE=T
}
\caption{Wilcoxon rank sum test, paired, one sided, significance level 5\%. Shown is the p-value for the hypothesis that  for a specific G-problem the full SACOBRA method at the final iteration is better than another solver M$*$. Significant improvements ($p \leq 5\%$) are marked as cells with dark blue color. Optimization methods: M1: SACOBRA$\backslash$rescale (SACOBRA without rescaling the input space), M2: SACOBRA$\backslash$RS (SACOBRA without random start), M3: SACOBRA$\backslash$aDRC, M4: SACOBRA$\backslash$aFF, M5: SACOBRA$\backslash$aCF, M6: COBRA ($\Xi=\Xi_{s}$), M7: COBRA ($\Xi=\Xi_{l}$).
%\SB{According to http://colorbrewer2.org/  I could not find any color coding with blue tonnage which are photocopy friendly for 4 colors, If we want to use definitely blue tonnage then we have to have only 3 classes}
%\WK{The statistical test should be made at the last iteration we see in the conv-curves. Is it so?}\SB{yes}
}%
\label{fig:wilcoxsb}%
\end{figure}

Fig.~\ref{fig:wilcoxsb} shows that each of these elements has its relevance for some of the G-problems: The full SACOBRA method is compared with other SACOBRA- or COBRA-variants M$*$ on 30 runs. SACOBRA is significantly better than each M$*$ at least for some G-problems (each column has a dark cell). And each G-problem benefits from one or more SACOBRA extensions (each row has a dark cell). The only exception from this rule is G11, but for a simple reason: 
%G02 is a high-dimensional and highly multimodal problem which is generally hard to solve by \textit{any} of the SACOBRA- or COBRA-variants (the reason is that surrogate models with a low number of points cannot capture enough detail of this complicated fitness function). 
G11 is an easy problem which is solved by \textit{all} SACOBRA variants in each run, so none is significantly better than the others.

\begin{figure}[tb]%
\begin{center}
\includegraphics[width=0.65\columnwidth]{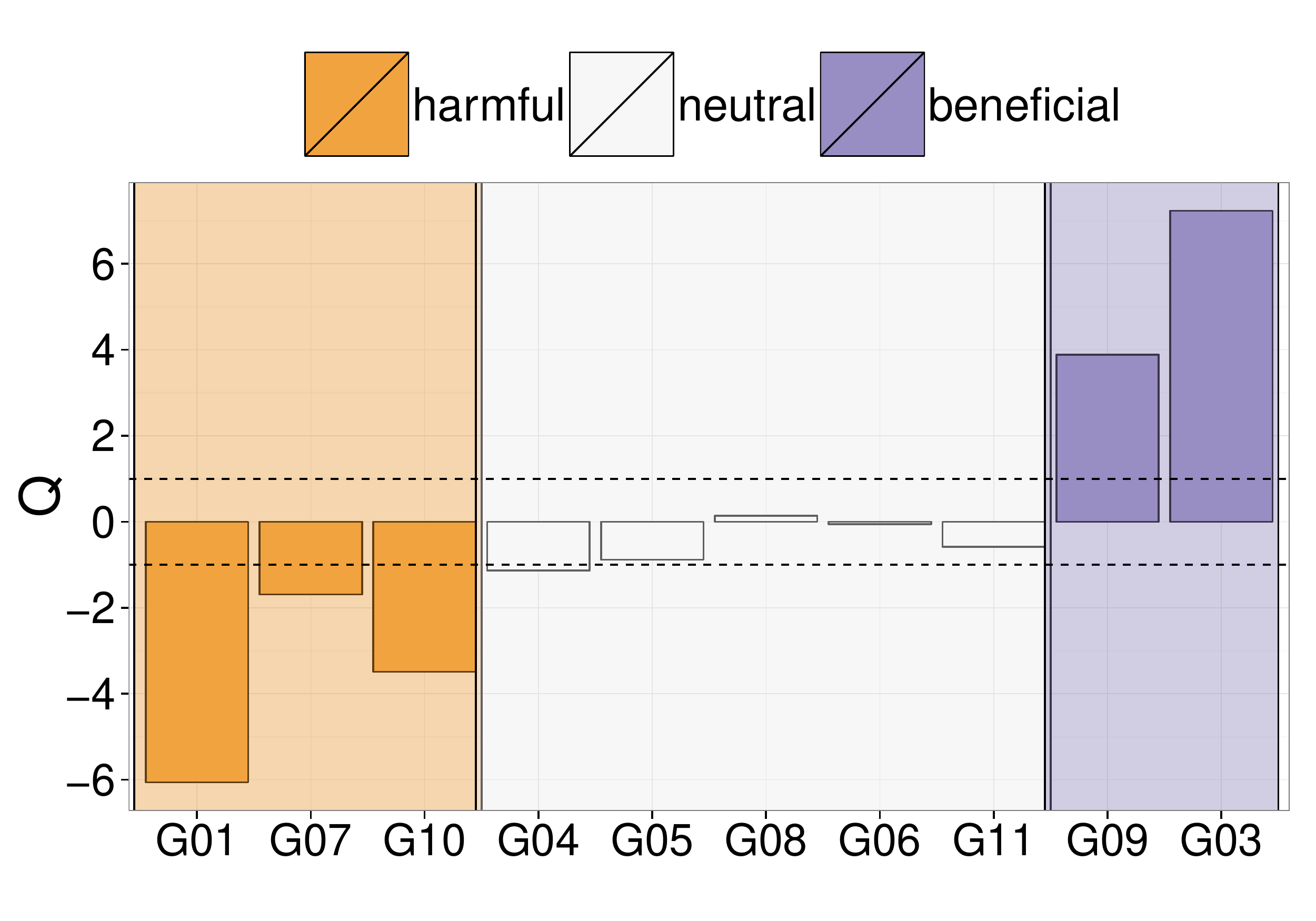}%
\end{center}
\caption{$Q$-value (Eq.~\eqref{eq:Qvalue}) at end of optimization for all G-problems. The G-problems are ordered along the x-axis according to the $R$-value defined in Eq.~\protect\eqref{eq:medErrRatio} which measures the impact of $plog$ on the optimization performance. Any threshold for $Q$ in $[-1,1]$ will clearly separate the harmful from the beneficial problems. This figure shows that the online available $Q$ is a good predictor of the impact of $plog$ on the overall optimization performance.
%\SB{These colors are supposed to be colorblind photocopy and print friendly colors, generated by:http://colorbrewer2.org/ }
%\WK{I would suggest to have only fig:qvalue in the final paper and integrate into it as a background to the coloured bars the rectangles of fig:peffectv1.} 
}%
\label{fig:qvalue}%
\end{figure}

%%%%%%%%%%%%%%%%%%%%%%%%%%%%%%%%%%%%%%%%%%%%%%%%%%%%%%%%%%%%%%%%%%%%%%%%%%%%%%%%%%%%%%%%%%%%%%%%%%%%
\subsection{Fitness Function Adjustment}
\label{Sec:fitfuncadjust}
%\SB{do we need this paragraph at all?} 
By comparing the convergence curves of G-functions we realized that applying the logarithmic transform is strictly harmful for three of the G-functions, significantly beneficial for two other problems, and with negligible effect on the other problems. Therefore, a careful selection should be done. Although we demonstrated in Sec.~\ref{Sec:plogPitfall} that steep functions can be better modeled after the logarithmic transformation, it is not trivial to define a correct threshold to classify steep functions. 
Also, there is no direct relation between steepness of the function and the effect of logarithmic transformation on optimization. 
We defined in Sec.~\ref{Sec:aFF} and Algorithm~\ref{alg:sacobra2Adjust}, function \textsc{AnalyzePlogEffect}, a measure called $Q$ in order to quantify online whether RBF models with and without $plog$ transformation are better or worse.

Here we test by experiments whether the $Q$-value does a good job.
Fig.~\ref{fig:qvalue} shows the $Q$-value for all G-problems. The G-problems are are ranked on the horizontal axis according to the impact of logarithmic transformation of the fitness function on the optimization outcome. This means that applying the $plog$-transformation has the worst effect for modeling the fitness of G01 and the best effect for G03. We measure the impact on optimization in the following way: For each G-problem we perform 30 runs with $plog$ inactive and with $plog$ active. We calculate the median of the final optimization error in both cases and take the ratio
\begin{equation}
		R = \frac{\med(E_{opt})}{\med(E_{opt}^{(plog)})}.
\label{eq:medErrRatio}
\end{equation}
Note that $R$ is usually not available in normal optimization mode.
If $R$ is \{~close to zero, close to 1, much larger than 1~\} then the effect of $plog$ on optimization performance is \{~harmful, neutral, beneficial~\}. It is a striking feature of Fig.~\ref{fig:qvalue} that the $Q$-ranks are very similar to the $R$-ranks.\footnote{The only notable difference, namely the switch in the order of G07 and G10, can be seen as an imperfection of measure $R$. Although G10 has rank 3 in $R$, it has weaker worst-case behavior than G07 because two G10 runs never produce a feasible solution if $plog$ is active.} This means that the beneficial or harmful effect of $plog$ is strongly  correlated with the RBF approximation error. 

Our experiments have shown that for all problems with $Q \in \left[-1,1\right]$ the optimization performance is only weakly influenced by the logarithmic transformation of the fitness function.  Therefore, in Step 19 of function \textsc{AdjustFitnessFunction} in Algorithm~\ref{alg:sacobra2Adjust}, any threshold in $\left[-1,1\right]$ will work. We choose the threshold 1, because it has the largest margin to the colored bars in Fig.~\ref{fig:qvalue}.

The G-problems for which $plog$ is beneficial are G03 and G09: These are according to Table~\ref{tab:GprobFeatures} the two problems with the largest fitness function range $FR$, thus strengthening our hypothesis from Sec.~\ref{Sec:plogPitfall}: For such functions a $plog$-transform should be used to get good RBF-models. The G-problems for which $plog$ is harmful are G01, G07, and G10: Looking at the analytical form of the objective function in those problems\footnote{The analytical form is available in the appendices of \protect\cite{runarsson2000stochastic} or \protect\cite{runarsson2005search}.} we can see that these are the only three functions being of quadratic type (Table~\ref{tab:GprobFeatures}) \textit{and} having no mixed quadratic terms. Those functions can be fitted perfectly by the polynomial tail (Eq.~\eqref{eq:squares}) in SACOBRA, if $plog$ is \textit{in}active. With $plog$ they become nonlinear and a more complicated approximation by the radial basis functions is needed. This results in a larger approximation error. 

%---------------------------------------------------------------------------------------------------
\input{longTab}					%% Table \ref{tab:bigTab} 
%---------------------------------------------------------------------------------------------------

%%%%%%%%%%%%%%%%%%%%%%%%%%%%%%%%%%%%%%%%%%%%%%%%%%%%%%%%%%%%%%%%%%%%%%%%%%%%%%%%%%%%%%%%%%%%%%%%%%%%
\subsection{Comparison with other optimizers}
\label{Sec:CompareOptim}

Table~\ref{tab:bigTab} shows the comparison with different state-of-the-art optimizers on the G-problem suite. While ISRES (Improved Stochastic Ranking \cite{runarsson2005search}) and DE (Differential Evolution \cite{brest2006saDE}) are the best optimizers in terms of solution quality, they require the highest number of function evaluations as well. SACOBRA has on most G-problems (except G02) the same solution quality, only G09 and G10 are very slightly worse. At the same time SACOBRA requires only a small fraction of function evaluations (fe): roughly 1/1000 as compared to ISRES and RGA and 1/300 as compared to DE (row \textit{average fe} in Table~\ref{tab:bigTab}).

G02 is marked in red cell color in Table~\ref{tab:bigTab} because it is not solved to the same level of accuracy by most of the optimizers. ISRES and RGA (Repair GA \cite{chootinan2006constraint}) get close, but only after more than 300\,000 fe. DE performs even better on G02, but requires more than 200\,000 fe as well. SACOBRA and COBRA cannot solve G02.

The results in column SACOBRA, DE and COBYLA are from our own calculation in \texttt{R} while the results in column COBRA, ISRES and RGA were taken from the papers cited. In two cases (red italic numbers in Table~\ref{tab:bigTab}) the reported solution is better than the true optimum, possibly due to a slight infeasibility. This is explicitly stated in the case of ISRES~\cite[p. 288]{runarsson2000stochastic}, because the equality constraint $h(x)=0$ of G03 is transformed into an approximate inequality $|h(x)| \leq \epsilon $ with $\epsilon=0.0001$.

COBRA \cite{regis2014constrained} comes close to SACOBRA in terms of efficiency (function evaluations), but it has to be noted that \cite{regis2014constrained} does not present results for all G-problems (G01 and G11 are missing and G02 results are for 10 dimensions, but the commonly studied version of G02 has 20 dimensions).
%(G1, G2, and G11 are missing)}. 
Furthermore, for many G-problems (G03, G06, G07, G09, G10) a manual transformations of the original fitness function or the constraint functions was done in \cite{regis2014constrained} prior to optimization.  SACOBRA starts without such transformations and proposes instead self-adjusting mechanisms to find suitable transformations (after the initialization phase or on-line). 

COBYLA often produces slightly infeasible solutions, these are the numbers in brackets. If such infeasible runs occur, the median was only taken over the remaining feasible runs, which is in principle too optimistic in favor of COBYLA. 

%---------------------------------------------------------------------------------------------------------------
%------comparing SACOBRA with DEoptim and COBYLA
%---performance and data profiles
% this plot is generated by compareCobylaDeoptim.R script (both figures from this script,
% with plotPP(1,50,  ...) for the  left plot DP5compareCobylaDeoptim0.pdf
% and  plotPP(1,1000,...) for the right plot DP5compareCobylaDeoptim1.pdf
%
% methods<-c("SACOBRA50","SACOBRA9","COBYLA","DEoptim")
%  for        SACOBRA     COBRA-R    COBYLA   DE
%
\begin{figure}[tb]%
\begin{minipage}[b]{0.5\textwidth}
	\includegraphics[width=\columnwidth]{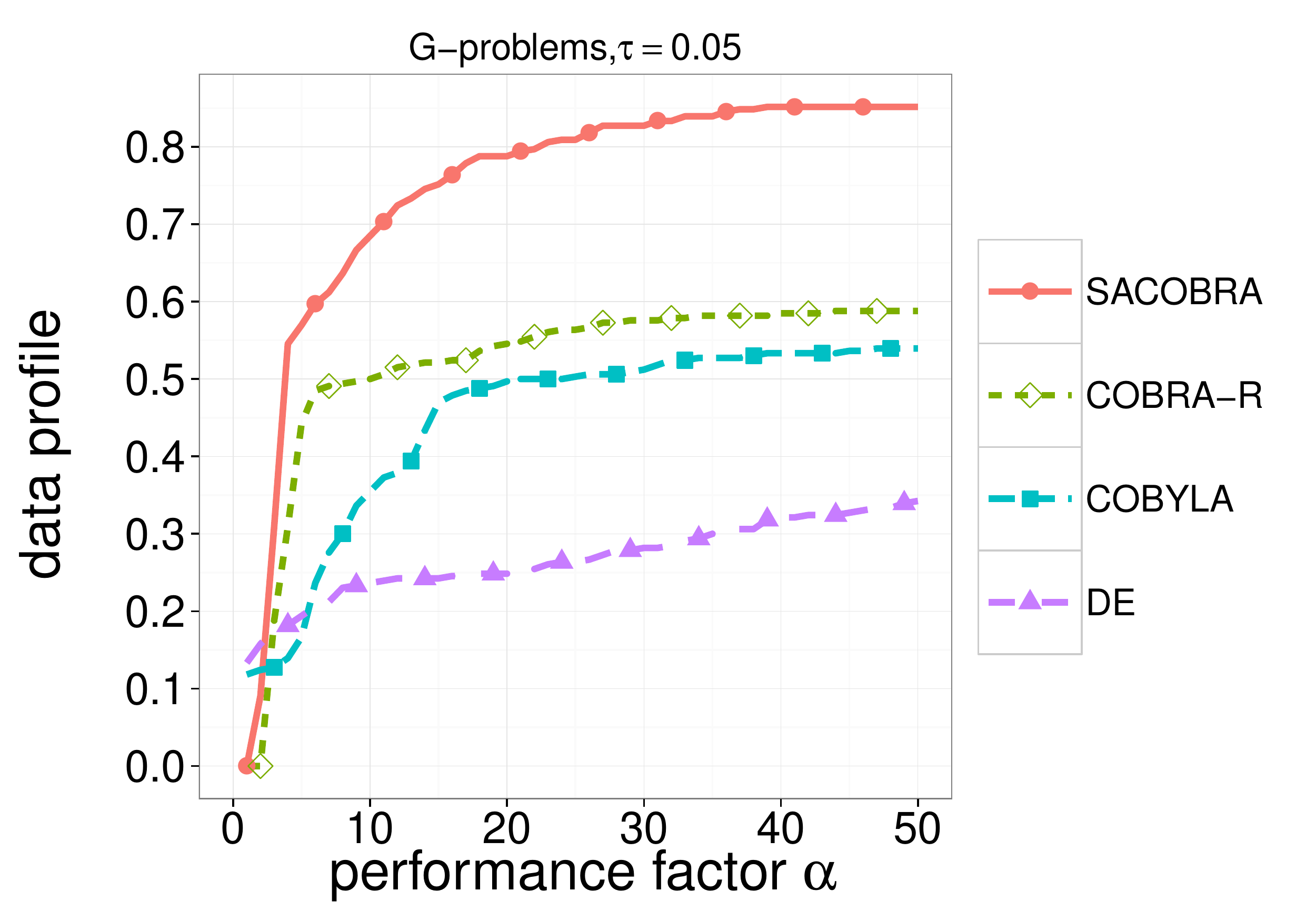}%
\end{minipage}
\hfill
	\begin{minipage}[b]{0.5\textwidth}
	\includegraphics[width=\columnwidth]{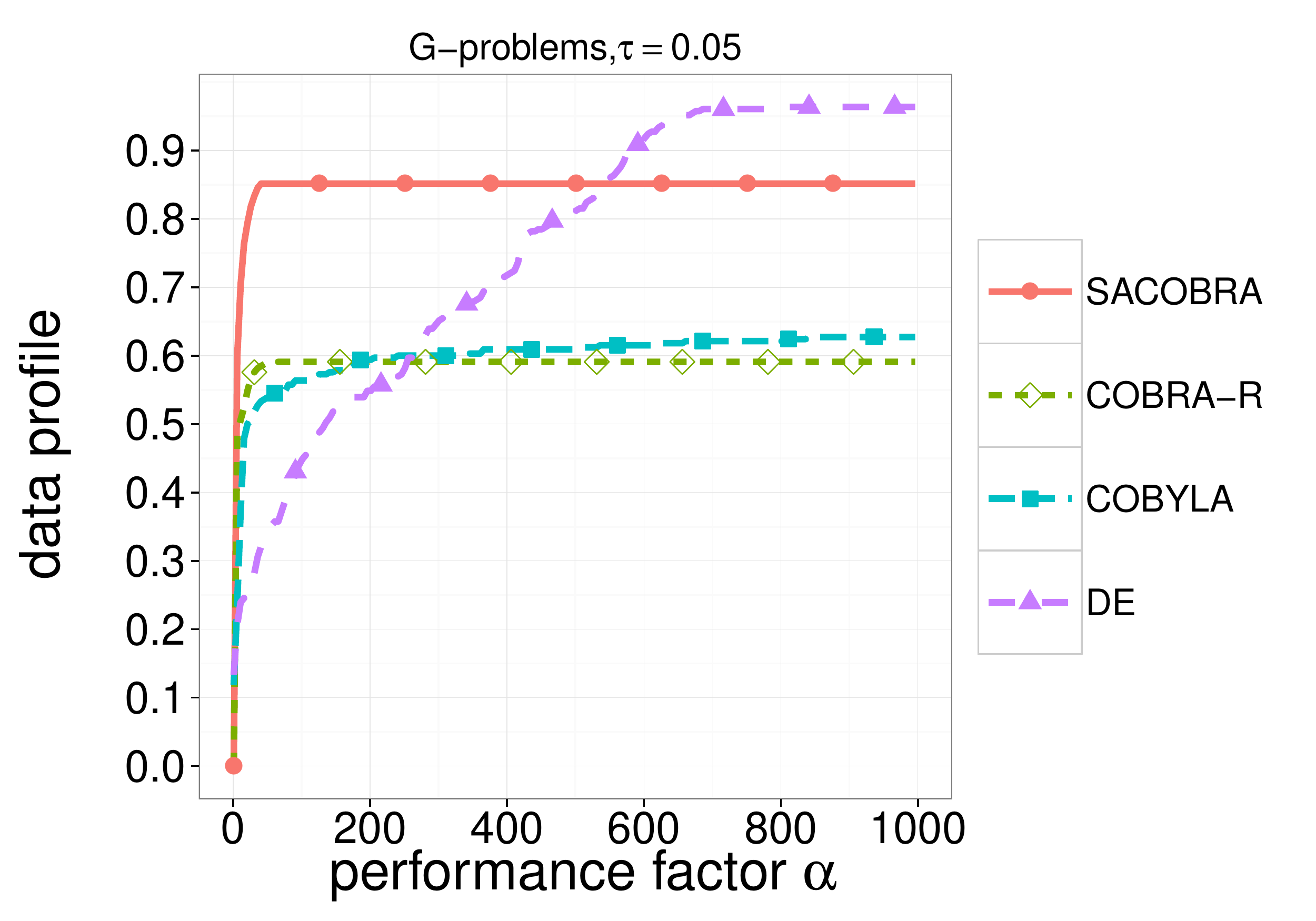}%
	\end{minipage}
%\end{tabular}
\caption{Comparing the performance of the algorithms SACOBRA, COBRA-R (with rescale), Differential Evolution (DE), and COBYLA on optimizing all G-problems G01-G11  (30 runs with different initial random populations).
%\WK{replace \textit{tolerance} by $\tau$ in the plot title}
%\WK{Since we do no longer distinguish between SACOBRA and SACOBRA2, can you prepare modified plots where the current SACOBRA is removed and the current SACOBRA2 is renamed to SACOBRA?}
%\WK{Can you make both plots identical in height?} 
%\SB{Changed}
%\WK{I would rename in the legend DEoptim to DE and COBRA-R (with rescale) simply to COBRA-R, so that the legend becomes a bit less wide and the plot area can be a bit bigger.}
%\WK{x label: there should be a space between "factor" and "$\alpha$"}
%
%\WK{A question concerning DE and the right figure: $\alpha=1000$ would mean only 21.000 fe's for G02-20d. Is this true, that DE can solve G02-20d after only 21.000 fe's? This would be much faster than ISRES! Or is it only that DE is better than the others after 21.000 fe's, but has not yet solved G02-20d to good accuracy?}\SB{What we are seeing is that DE is better than the other algorithms after $1000* (d+1)$ but still does not reach the 100\%.  The average required fe's to solve G02-20d is 226994 which is quite a large number but still 100 thousand iterations less than what ISRES needs.}
}%
\label{fig:sacobraDeCobyla}%
\end{figure}
%\end{center}
%---------------------------------------------------------------------------------------------------------------

\begin{table}[tb]%
\caption{Number of infeasible runs among 330 runs returned by each method on the G-problem benchmark. A run is infeasible if the final best solution is infeasible.
}
\label{tab:infeasRun}
\centering
\begin{tabular}{lcc}
\\\hline
method                     &  %\parbox[][][c]{2.5cm}{
															infeasible runs %}
															& functions\\\hline
\hline
SACOBRA                    & 0                & --\\
SACOBRA$\backslash$ rescale& 4                & G05\\
SACOBRA$\backslash$ RS     & 13               & G03, G05, G07,G09,G10\\
SACOBRA$\backslash$ aDRC   & 0                & --\\
SACOBRA$\backslash$ aFF    & 1                & G10\\
SACOBRA$\backslash$ aCF    & 0                & --\\
COBRA-R(no rescale)      	 & 37               & G03,G05,G07,G09,G10\\
COBRA-R(rescale)         	 & 23               & G05,G07,G09,G10\\
COBYLA										 & 48								& G02,G03,G05,G06,G07,G09,G10\\ 
DE												 &	0								& -- \\\hline
\end{tabular}
\end{table}\textbf{}

\begin{figure}[tb]
\centering
\includegraphics[width=0.8\columnwidth]{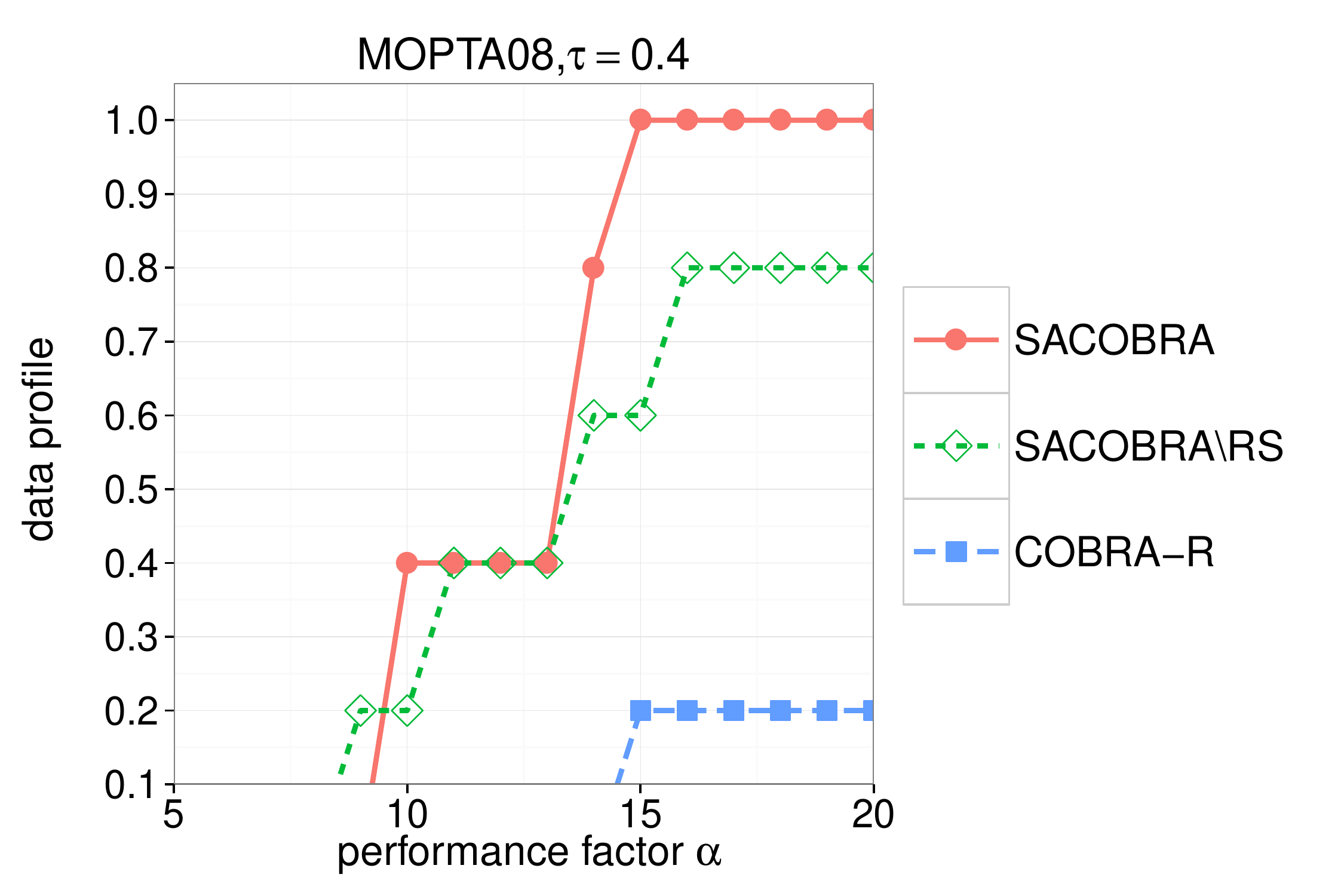}%
\caption{Data profile for MOPTA08: Same as Fig.~\ref{fig:SACOBRAcompfinal} but with 10 runs on MOPTA08 with different initial designs. % (instead of 30 runs on all 11 G-problems).
The curves for SACOBRA without rescale, aDRC, aFF, or aCF are identical to full SACOBRA, since in the case of MOPTA08 the objective function and the constraints are already normalized. 
%Analyzing the impact of different elements of SACOBRA on MOPTA08. Data profile of SACOBRA, SACOBRA$\backslash$rescale (SACOBRA without rescaling the input space), and other \glqq $\backslash$\grqq-algorithms are  with a similar meaning. COBRA-R is the old version of COBRA~\cite{Koch14a}, i.~e. SACOBRA with all adjustment extensions switched off. These algorithms are performed on 10 independently initialized runs. 
%\SB{aDRC, aFF, aCF are identical}
%\WK{Please reverse the legend for the final plot (SACOBRA at the top)}
}%
\label{fig:moptasacobra}%
\end{figure}

Fig.~\ref{fig:sacobraDeCobyla} shows the comparison of SACOBRA and COBRA-R with other well-known constraint optimization solvers available in R, namely DE\footnote{\texttt{R}-package \texttt{DEoptimR}, available from 
\href{https://cran.r-project.org/web/packages/DEoptimR}{\url{https://cran.r-project.org/web/packages/DEoptimR}}}
 and COBYLA.\footnote{\texttt{R}-package \texttt{nloptr}, available from 
\href{https://cran.r-project.org/web/packages/nloptr}{\url{https://cran.r-project.org/web/packages/nloptr}}}
%{https://cran.r-project.org/web/packages/nloptr}. 
The right plot in Fig.~\ref{fig:sacobraDeCobyla}  
shows that DE achieves good results after many function evaluations, in accordance with Table~\ref{tab:bigTab}. But the left plot in Fig.~\ref{fig:sacobraDeCobyla} shows that DE is not really competitive if very tight bounds on the budget are set.

Tab.~\ref{tab:infeasRun} shows that SACOBRA greatly reduces the number of infeasible runs as compared to COBRA-R.  Most of the SACOBRA variants have less than 2\% infeasible runs whereas COBRA-R has 7-11\%. The full SACOBRA method has no infeasible runs at all. 

\subsection{MOPTA08}
Fig.~\ref{fig:moptasacobra} shows that we get good results with SACOBRA on the high-dimensional MOPTA08 problem ($d=124$) as well. 
 A problem is said to be \textit{solved} in the data profile of Fig.~\ref{fig:moptasacobra} if it is not more than $\tau=0.4$ away from the best value obtained in all runs by all algorithms. 

Table~\ref{tab:mopta} shows the results after 1000 iterations for Regis' recent trust-region based approach TRB~\cite{regis2015trust} and our algorithms. We can improve the already good mean best feasible results of 227.3 and 226.4 obtained with COBRA-R~\cite{Koch2015a} and TRB~\cite{regis2015trust}, resp., to 223.3 with SACOBRA. The reason that SACOBRA$\backslash$RS is slightly better than COBRA-R~\cite{Koch2015a} is that SACOBRA uses an improved DRC.
%\WK{SB, can you confirm these numbers}\SB{done}

\begin{center}
\begin{table*}[tb]%
\caption{Comparing different algorithms on optimizing MOPTA08 after 1000 function evaluations. }
\centering
\begin{tabular}{lrrrr}
Algorithm 								&best &median &mean	&worst \\\hline
COBRA-R~\cite{Koch2015a}	&226.3&227.0 &227.3&229.5\\
TRB~\cite{regis2015trust}	&225.5&226.2&226.4&227.4\\
SACOBRA$\backslash$RS			&\textbf{222.4}&\textbf{223.1} &223.6 &224.8\\
SACOBRA										&223.0&223.3&\textbf{223.3}&\textbf{223.8}\\ \hline
\end{tabular}
\label{tab:mopta}
\end{table*}
\end{center}
%%%%%%%%%%%%%%%%%%%%%%%%%%%%%%%%%%%%%%%%%%%%%%%%%%%%%%%%%%%%%%%%%%%%%%%%%%%%%%%%%%%%%%%%%%%%%%%%%%%%
\section{Discussion}
%%%%%%%%%%%%%%%%%%%%%%%%%%%%%%%%%%%%%%%%%%%%%%%%%%%%%%%%%%%%%%%%%%%%%%%%%%%%%%%%%%%%%%%%%%%%%%%%%%%%
\label{Sec:Discussion}

\subsection{SACOBRA and surrogate modeling}
SACOBRA is an algorithm capable of self-adjusting its parameters to a wide-ranging set of problems in constraint optimization. We analyzed the different elements of SACOBRA and their importance for efficient optimization on the G-problem benchmark. It turned out that the two most important elements are rescaling (especially in the early phase of optimization)  and automatic fitness function adjustment (aFF, especially in the later phase of optimization). Exclusion of either one of these two elements led to the largest performance drop in Fig.~\ref{fig:SACOBRAcompfinal} compared to the full SACOBRA algorithm.

We may step back for a moment and ask why these two elements are important. Both of them are directly related to accurate RBF modeling, as our analysis in Sec.~\ref{Sec:rescale} %and Sec.~\ref{Sec:plogPitfall} 
has shown. If we do not rescale, then the RBF model for a problem like G10 will have large approximation errors due to numeric instabilities.  If we do not perform the $plog$-transformation in problems like G03 with a very large fitness range $FR$ (Tab.~\ref{tab:GprobFeatures}) and thus very steep regions, then such problems cannot be solved. This can be attributed to large RBF approximation errors as well. 

%%%% The following paragraph is probably no longer true with the new results from SACOBRA5x. And we find no clear difference for "approximation 
%%%% error at the true solution" in the cases "RS" and "no RS".
%
%The third element found to be important is RS. If this element is important it means that the optimization would otherwise get stuck in a suboptimal local minimum. This can happen if the RBF model has a local minimum in the wrong region and either does not model the true minimum region appropriately or the optimizer cannot find a way to the true minimum.
%\WK{SB, can you have a look at the cases where we have bad runs in the case of "no RS": Is it true that in those cases the approximation error of the RBF model at the true optimum is constantly bad? If so, we could argue that RS solves an RBF-modelling-problem as well.} \SB{difficult to say, please see Figure~\ref{fig:RSanalyse}} 

%\WK{moved this paragraph from Conclusion to Discussion, it is too detailed for a conclusion}
We diagnosed that the quality of the surrogate models \replacedOK{is in relationship}{ has a relation} with the correct choice of the DRC parameter, which controls the step size in each iteration. It is more desirable to choose a set of smaller step sizes for functions with steep slopes. An automatic adjustment step in SACOBRA can identify steep functions after a few function evaluations and decide whether to use a large DRC or a small one.

%\WK{moved this paragraph from Conclusion to Discussion, it is too detailed for a conclusion}
For the G-problem suite the constraint functions vary in number, type and range. Our experiments showed that handling all constraints can be challenging, especially when the constraint functions have widely different ranges. For that reason, we considered an automatic adjustment approach to normalize all the constraints by using the information gained about the constraints after the evaluation of the initial population. The SACOBRA algorithm also benefits from using a random start mechanism to avoid getting stuck in a local optimum of the fitness surrogate. 

%\WK{TODO: review self-adjusting elements SACOBRA, comparison with results of other self-adjusting algos, e.g. Farmani~\cite{farmani2003saFitness}.} 
%\WK{not needed anymore}

\subsection{Limitations of SACOBRA}
\label{sec:limitation}
\subsubsection{Highly multimodal functions}
Surrogate models like RBF are a great thing for efficient optimization and probably the only way to solve constrained optimization problems in less than 500 iterations. But a current border for surrogate modeling are highly multimodal functions. G02 is such a function, it has a large number of local minima. Those functions have usually large first and higher order derivatives. If a surrogate model interpolates isolated points of such a function, it tends to overshoot in other parts of the function. 
To the best of our knowledge, highly multimodal problems cannot be solved so far by surrogate models, at least not for higher dimensions with high accuracy. This is also true for SACOBRA.  
%\WK{Tenne and Armfield~\cite{tenne2008} present an approach where e.g. Rastrigin 30d is solved in less than 200 fe. But they do not have precision, the error is in the order of 50 or more (!)}
Usually the RBF model has a good approximation only in the region of one of the local minima and a bad approximation in the rest of the search space.  Further research on highly multimodal function approximation is required to solve this problem. 

\subsubsection{Equality constraints} 
The current approach in COBRA (and in SACOBRA as well) can only handle inequality constraints. The reason is that equality constraints do not work together with the uncertainty mechanism of Sec.~\ref{Sec:Cobra}. A reformulation of an equality constraint $h(x)=0$ as inequality 
$|h(x)| \leq 0.0001$ as in~\cite{runarsson2000stochastic}  is not well-suited for COBRA and for RBF modeling. We used in this work the same approach as Regis~\cite{regis2014constrained} and replaced each equality operator with an inequality operator of the appropriate direction. It has to be noted however, that such an approach contradicts a true black-box handling of constraints and that it is -- although being viable for the problems G01-G11 -- not viable for more complicated objective functions having their minima on both sides of equality constraints. In a forthcoming paper~\cite{bagheri2016equ} we will address this problem separately.

%%%%%%%%%%%%%%%%%%%%%%%%%%%%%%%%%%%%%%%%%%%%%%%%%%%%%%%%%%%%%%%%%%%%%%%%%%%%%%%%%%%%%%%%%%%%%%%%%%%%
\section{Conclusion}
%%%%%%%%%%%%%%%%%%%%%%%%%%%%%%%%%%%%%%%%%%%%%%%%%%%%%%%%%%%%%%%%%%%%%%%%%%%%%%%%%%%%%%%%%%%%%%%%%%%%
\label{Sec:Conclusion}

We summarize our discussion by stating that a good understanding of the capabilities and limitations of RBF surrogate models -- which is not often undertaken in the surrogate literature we are aware of -- is an important prerequisite for efficient and effective constrained optimization. 

The analysis of the errors and problems occurring initially for some G-problems in the COBRA algorithm have given us a better understanding of RBF models and led to the development of the enhancing elements in SACOBRA.
\replacedOK{By studying a widely varying set of problems we observed certain challenges when modeling very steep or relatively flat functions with RBF. This can result in large approximation errors.}{
Using a widely varying set of problems for benchmarking RBF surrogates, assisted us to observe the challenges for modeling steep functions due to the numerical instabilities.}  SACOBRA tackles this problem by making use of \replacedOK{a conditional $plog$-transform for the objective function}{the online fitness function $plog$-transform}. \addOK{We proposed a new online mechanism to let SACOBRA decide automatically when to use $plog$ and when not.}

Numerical issues to train RBF models can also occur in the case of a very large input space. A simple solution to this problem is to rescale the input space. Although many other optimizers recommend to rescale the input, this work has shown the reason behind it and the importance of it by evidence.  Therefore, we can answer our first research question \textbf{(H1)} positively\addOK{: Numerical instabilities can occur in RBF modeling, but it is possible to avoid them with the proper function transformations and search space adjustments}.

SACOBRA benefits from all \replacedOK{its extension elements introduced in Sec.~\ref{Sec:Sacobra}}{ mentioned elements}. Each element boosts up the optimization performance on a subset of \addOK{all} problems without harming the optimization process on the other ones. As a result, the overall optimization performance on the whole set of problems is improved by 50\% as compared to COBRA (with a fixed parameter set). About 90\% of the tested problems can be solved efficiently by SACOBRA\addOK{ (Fig.~\ref{fig:SACOBRAcompfinal})}. %Accordingly, \replacedOK{our}{the} answer \replacedOK{to question}{ of} 
The answer to \textbf{(H2)} is: SACOBRA is capable to cope with many diverse challenges in constraint optimization. \addOK{It is the main contribution of this paper to propose with SACOBRA the first surrogate-assisted constrained optimizer which solves efficiently the G-problem benchmark and requires \textbf{no parameter tuning} or manual function transformations.} Finally, let us provide a result to \textbf{(H3)}:  SACOBRA requires \textbf{less than 500 function evaluations} to solve 10 out of 11 G-problems (exception: G02) with similar accuracy as \replacedOK{other state-of-the-art algorithms}{ the other state of the arts}. Those other algorithms often need between 300 and 1000 times more function evaluations.

%\replaced[remark=\WK{In principle, all what is said in the replaced paragraph was already said in Sec.~\ref{sec:limitation}. Therefore I tried to formulate it much shorter.}]
%{
Our future research will be devoted to overcome the current limitations of SACOBRA mentioned in Sec.~\ref{sec:limitation}. These are: (a) highly multimodal functions like G02 and (b) equality constraints.
%}{A demanding condition in efficient constraint optimization is to model highly multimodal objective or constraint functions with only a few points. 
%At present, SACOBRA \replacedOK{has its limitations in dealing with highly multimodal functions}{ cannot address such problems} in an efficient way. 
%Our results for \replacedOK{G02 show}{ the G02 problem shows} the poor performance of SACOBRA when the objective function is highly multimodal. 
%Our future investigation will focus on tackling this class of problems. On the other hand, SACOBRA is limited to address only inequality constraints. Further work to eliminate this limitation directs us to develop an optimization framework which can cover a more diverse set of problems.}

\section*{Acknowledgements}
\vspace{-0.9cm}
\begin{tabular}{p{0.7\columnwidth}r}
\vspace{-2cm}
This work has been supported by the Bundesministerium f\"{u}r Wirtschaft (BMWi) under the ZIM grant MONREP (AiF FKZ KF3145102, Zentrales Innovationsprogramm Mittelstand). 
&
\includegraphics[width=0.25\columnwidth]{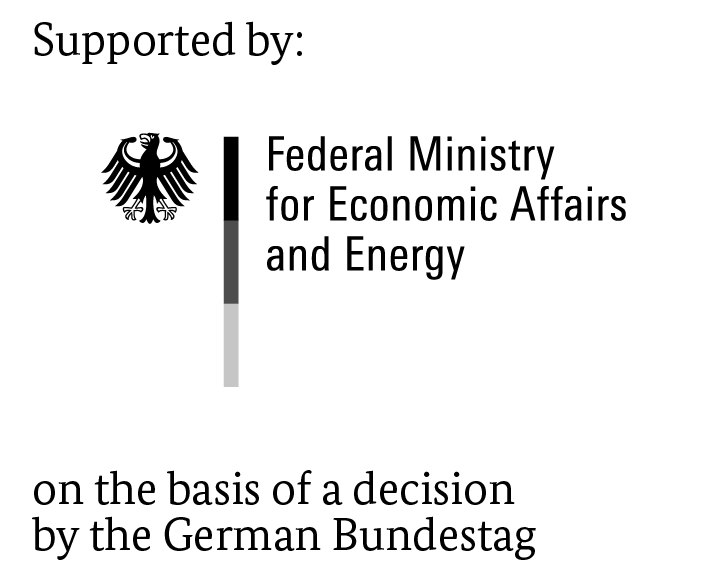}
\end{tabular}

%\bibliography{constraints}
%\bibliographystyle{abbrv}

\input{arXiv-SACOBRA15-bibitems}
\section*{Vitae}

\noindent
{\includegraphics[width=1in,height=1.25in,keepaspectratio]{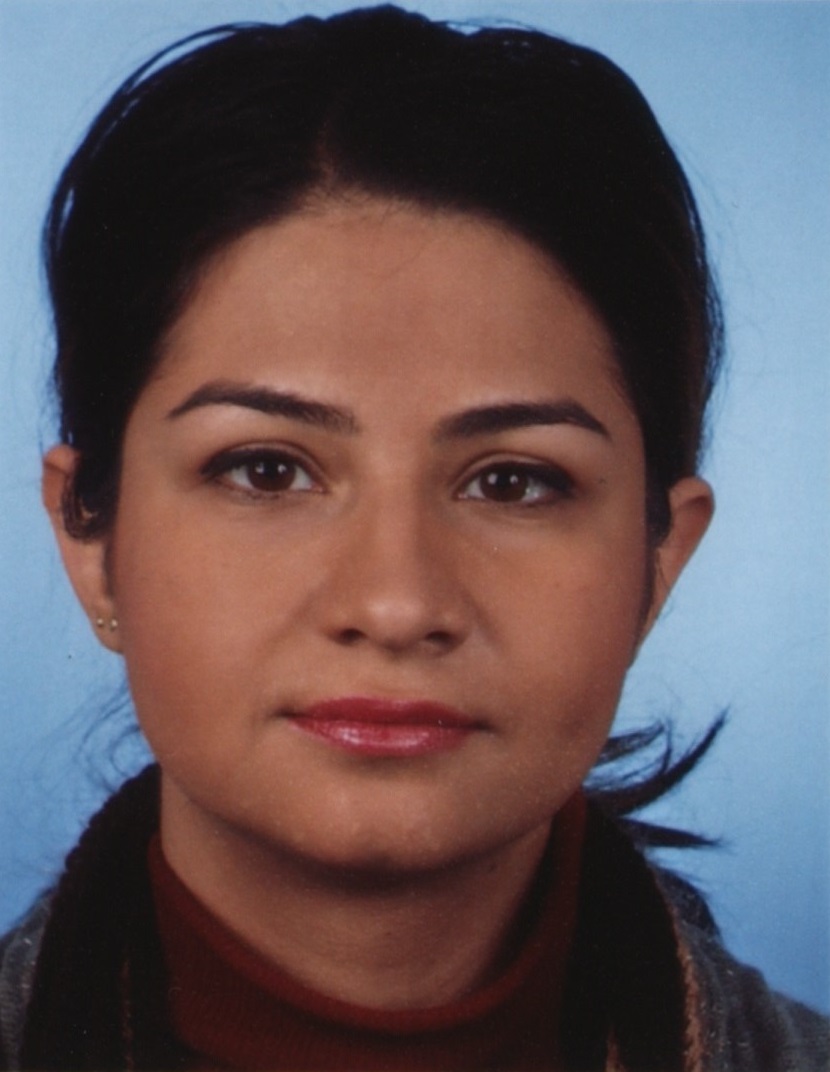}}
{\bf Samineh Bagheri}
received her B.Sc. degree in electrical engineering specialized in electronics from Shahid Beheshti University, Tehran, Iran in 2011. She received her M.Sc. in Industrial Automation \& IT from Cologne University of Applied Sciences, where she is currently research assistant and Ph.D. student in cooperation with Leiden University.
\\Her research interests are machine learning, evolutionary computation and constrained and multiobjective optimization tasks.
\vspace{1cm}

\noindent
\includegraphics[width=1in,height=1.25in,keepaspectratio]{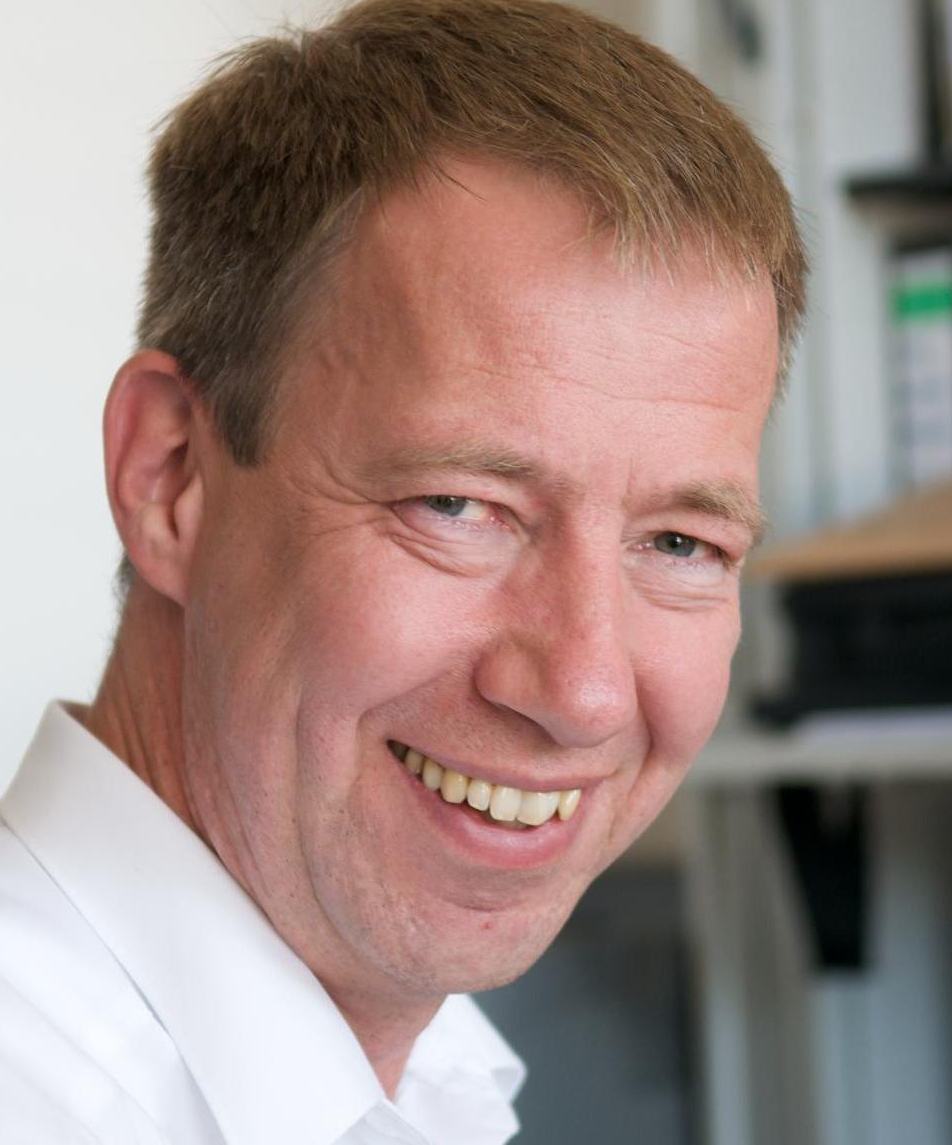}
{\bf Wolfgang Konen}
is Professor of Computer Science and Mathematics at Cologne University of Applied Sciences, Germany.
He received his Diploma in physics and his
Ph.D. degree in theoretical physics from the
University of Mainz, Germany, in 1987 and 1990, resp. He worked in the area of neuroinformatics and computer vision 
at Ruhr-University Bochum, Germany, and in several companies. % from 1990--2002.
He is 
founding member of the Research Centers Computational Intelligence, Optimization \& Data Mining (\url{http://www.gociop.de}) and CIplus (\url{http://ciplus-research.de}). He co-authored more than 100 papers and 
his research interests include, but are not limited to: efficient optimization, neuroevolution, machine learning, data mining, and computer vision.
\vspace{1cm}

\noindent
{\includegraphics[width=1in,height=1.25in,clip,keepaspectratio]{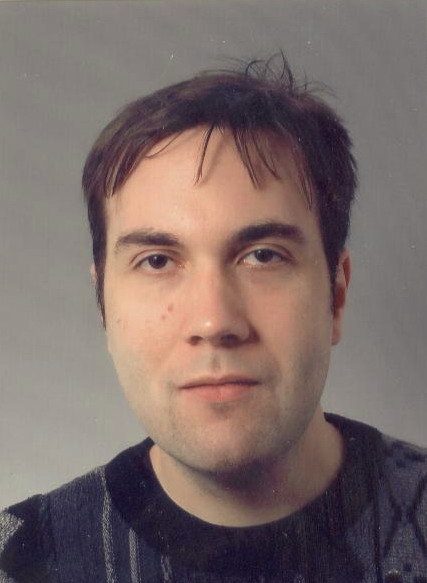}}
{\bf Dr. Michael Emmerich} is Assistant Professor at LIACS, Leiden University, and leader of the Multicriteria Optimization and Decision Analysis research group. He received his doctorate in 2005 from Dortmund University (H.-P. Schwefel, promoter) and worked as a researcher at ICD e.V. (Germany), IST Lisbon, University of the Algarve (Portugal), ACCESS Material Science e.V. (Germany), and the FOM/AMOLF institute (Netherlands). He is known for pioneering work on model-assisted and indicator-based multiobjective optimization, and has co-authored more than 100  papers in machine learning, multicriteria optimization and surrogate-assisted optimization and its applications in chemoinformatics and engineering optimization.
\vspace{1cm}

\noindent
{\includegraphics[width=1in,height=1.25in,clip,keepaspectratio]{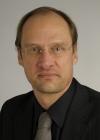}}
{\bf Thomas B\"ack} is head of the Natural Computing Research Group at the Leiden Institute of Advanced Computer Science (LIACS). He received his PhD in Computer Science from Dortmund University, Germany, in 1994. He has been Associate Professor of Computer Science at Leiden University since 1996 and full Professor for Natural Computing since 2002. 
Thomas B\"ack has more than 150 publications on natural computing technologies. His main research interests are theory and applications of evolutionary algorithms (adaptive optimization methods gleaned from the model of organic evolution), cellular automata, data-driven modeling and applications of those methods in medicinal chemistry, pharmacology, and engineering. 

\end{document}

%% file: cobraFlc.tex
\tikzstyle{startstop} = [circle, minimum size=0.7cm,text centered, draw=black!50, fill=black]
\tikzstyle{io} = [trapezium, trapezium left angle=70, trapezium right angle=110, text width=3cm, minimum height=1cm, text centered, draw=black, fill=white]
\tikzstyle{process} = [rectangle, minimum width=3cm, minimum height=1cm, text centered, text width=3cm, draw=black, rounded corners, fill=white]
\tikzstyle{processN} = [rectangle, minimum width=3cm, minimum height=1cm, text centered, text width=3cm, draw=black, rounded corners, fill=gray!50!white]
\tikzstyle{processS} = [rectangle, minimum width=3cm, minimum height=1cm, text centered, text width=3cm, draw=black, rounded corners, fill=orange]
\tikzstyle{decision1} = [diamond, aspect=1.8, inner sep=-2pt, minimum width=1cm, minimum height=1cm, text centered, text width=2cm, draw=black, fill=white]		% inner sep=-2pt brings the text closer to the border of the diamond
\tikzstyle{decision2} = [diamond, aspect=1.8, inner sep=-1pt, minimum width=1cm, minimum height=1cm, text centered, text width=2cm, draw=black, fill=white]		
\tikzstyle{arrow} = [->,very thick,>=stealth]

\begin{tikzpicture}[node distance=2cm and 4cm,on grid=true]
%\draw[step=1cm,gray,very thin] (0,0) grid (10,-10);

	\node (start) [startstop] {};
	\node (gen1) [process, right =of start] {Generate \& evaluate initial design}; 
	\node (repai) [process, below of=start] {Run repair heuristic};
	\node (deci1) [decision1, below of=repai] {Solution repaired or feasible?};
	\node (eva1) [process, below of =gen1, yshift=-2.0cm] {Evaluate solution on real functions};
	\node (upda) [process, below =of eva1] {Update the best solution};
	\node (deci2) [decision2, below of =upda] {Budget exhausted?};         
	\node (add1) [process, above =of eva1] {Add solution to the population}; 	 
	\node (opt1) [process, right =of eva1] {Run optimization on surrogates};   
	\node (fit1) [process, right =of deci2] {Fit RBF surrogates of objective and constraints};   
	\node (stop) [startstop, left =of deci2] {};			

	\draw [arrow] (start) -- (gen1);
	\draw [arrow] (gen1) --++(+171pt,0) |- (fit1);
	\draw [arrow] (add1) -- (eva1);
	\draw [arrow] (repai) -- (add1);
	\draw [arrow] (eva1) -- (deci1);
	\draw [arrow] (deci1) |- node[anchor=south,xshift=-0.38cm,yshift=0.3cm] {Yes} (upda);
	\draw [arrow] (deci1) -- node[anchor=north,xshift=-0.32cm,yshift=0.2cm] {No} (repai);
	\draw [arrow] (upda) -- (deci2);
	\draw [arrow] (deci2) -- node[anchor=west,yshift=-0.35cm,xshift=0.2cm] {Yes} (stop);
	\draw [arrow] (deci2) -- node[anchor=east,yshift=-0.30cm,xshift=0.2cm] {No} (fit1);
	\draw [arrow] (fit1) -- (opt1);
	\draw [arrow] (opt1) |- (add1);

\end{tikzpicture}

%% file: sacobra2Flc.tex
\tikzstyle{startstop} = [circle, minimum size=0.7cm,text centered, draw=black!50, fill=black]
\tikzstyle{io} = [trapezium, trapezium left angle=70, trapezium right angle=110, text width=3cm, minimum height=1cm, text centered, draw=black, fill=white]
\tikzstyle{process} = [rectangle, minimum width=3cm, minimum height=1cm, text centered, text width=3cm, draw=black, rounded corners, fill=white]
\tikzstyle{processN} = [rectangle, minimum width=3cm, minimum height=1cm, text centered, text width=3cm, draw=black, rounded corners, fill=gray!50!white]
\tikzstyle{processS} = [rectangle, minimum width=3cm, minimum height=1cm, text centered, text width=3cm, draw=black, rounded corners, fill=orange]
\tikzstyle{decision1} = [diamond, aspect=1.8, inner sep=-2pt, minimum width=1cm, minimum height=1cm, text centered, text width=2cm, draw=black, fill=white]		% inner sep=-2pt brings the text closer to the border of the diamond
\tikzstyle{decision2} = [diamond, aspect=1.8, inner sep=-1pt, minimum width=1cm, minimum height=1cm, text centered, text width=2cm, draw=black, fill=white]		
\tikzstyle{arrow} = [->,very thick,>=stealth]

\begin{tikzpicture}[node distance=2cm and 4cm,on grid=true]
%\draw[step=1cm,gray,very thin] (0,0) grid (10,-10);

%	\node (start) [startstop] {};
%	\node (resc) [processN, right =of start] {Rescale input space}; 
	\node (resc) [processN] {Rescale input space}; 
	\node (start) [startstop, left =3.2cm of resc] {};
	\node (gen1) [process, right =of resc] {Generate \& evaluate initial design}; 
	\node (sacobra1) [processN, right =of gen1] {Adjust constraint function(s)};
	\node (sacobra2) [processN, below of=sacobra1] {Adjust DRC};
	\node (repai) [process, below of=resc] {Run repair heuristic};
	\node (deci1) [decision1, below of=repai] {Solution repaired or feasible?};
	\node (upda) [process, below =of deci1] {Update the best solution};
	%\node (tr) [processS,below =of upda]{Refine the best solution};
	\node (deci2) [decision2, below of =upda] {Budget exhausted?};
	\node (fit) [processN,right =of deci2]{Online adjustment of fitness function};
	\node (stop) [startstop, left =3.2cm of deci2] {};	
	\node (fit1) [process, right =of fit] {Fit RBF surrogates of objective and constraints};   
  \node (rs) [processN, above=of fit1] {Select start point (RS)};
	\node (opt1) [process, above =of rs] {Run optimization on surrogates};   
	\node (add1) [process, right =of repai] {Add solution to the population}; 	 
	\node (eva1) [process, below of =add1] {Evaluate new point on real functions};
	
	\draw [arrow] (start) -- (resc);
	\draw [arrow] (resc) -- (gen1);
	\draw [arrow] (gen1)--(sacobra1);
	\draw [arrow] (sacobra1)--(sacobra2);
	\draw [arrow] (sacobra2) --++(+57pt,0) --++(0,-195pt) -| (fit);
	\draw [arrow] (deci2) -- node[anchor=south]{Yes} (stop);
	\draw [arrow] (deci2) -- node[anchor=south]{No} (fit);
	\draw [arrow] (fit) -- (fit1);
	\draw [arrow] (fit1) -- (rs);
	\draw [arrow] (rs) -- (opt1);
	%\draw [arrow] (opt1) -- (add1);
	\draw [arrow] (opt1.west) --++(-10pt,0)|- (add1.east);
	\draw [arrow] (add1)--(eva1);
	%\draw [arrow] (eva1.west) --++(-5pt,0)|- (deci1.east);
	\draw [arrow] (eva1)--(deci1);
	\draw [arrow] (deci1) -- node[anchor=west]{No}(repai);
	\draw [arrow] (deci1) -- node[anchor=west]{Yes}(upda);
	\draw [arrow] (upda) -- (deci2);
	%\draw [arrow] (tr) -- (deci2);
	\draw [arrow] (repai) -- (add1);

\end{tikzpicture}

%% file: ASOC-sacobra2Alg.tex
\begin{algorithm}%
\caption{SACOBRA. \textbf{Input:} %Deterministic 
Objective function $f$, set of constraint function(s) $\textbf{g}=\left( g_{1},\ldots,g_{m}\right) : 
[\vec{a},\vec{b}] \subset \mathbb{R}^d \rightarrow \mathbb{R}$ (see Eq.~\eqref{eq:minprob}),  
%which return a value at any point in the search space $\vec{x}\in [\vec{a},\vec{b}] \subset \mathbb{R}^d$. 
initial starting point $\vec{x}_{init}\in [\vec{a},\vec{b}] $, maximum evaluation budget $N_{max}$. 
\textbf{Output:} The best solution $\vec{x}_{best}$ found by the algorithm.}
\label{alg:sacobra2}
	\begin{algorithmic}[1] 
\Function{SACOBRA}{$f,\textbf{g},\vec{x}_{init},N_{max}$}
\State Rescale the input space to $\left[-1,1\right]^d$
\State Generate a random initial population: $P=\left\{x_{1}, x_{2},\cdots,x_{3 \cdot d}\right\}$% and evaluate $P$ on $f, \textbf{g}$
\State $(\widehat{FR},\widehat{GR}_i )$=\Call{AnalyseInitialPopulation}{$P,f,\textbf{g}$}
\State $\widetilde{\textbf{g}} \leftarrow$ \Call{AdjustConstraintFunctions}{$\widehat{GR}_i$, $\textbf{g}$}
\State $\Xi \leftarrow$ \Call{AdjustDRC}{$\widehat{FR}$}
\State $Q \leftarrow$ \Call{AnalysePlogEffect}{$f, P, \vec{x}_{init}$}
\State $\vec{x}_{best} \leftarrow \vec{x}_{init}$  %; $c \leftarrow 0$  %$F_{success} \leftarrow \textsc{True}$
\While{(budget not exhausted, $|P| < N_{max}$)}
\State $n \leftarrow |P|$
\State $\widetilde{f}() \leftarrow \left(\,Q > 1 ~?~ \ plog(f()) \,:\, f() \,\right)$          \Comment{see function $plog$ in Eq.~\eqref{eq:plog}}
\State Build surrogate models $\vec{s}\,^{(n)}$=$(s_{0}^{(n)},s_{1}^{(n)},\cdots,s_{m}^{(n)})$ for $(\widetilde{f}, \widetilde{g}_{1},\cdots,\widetilde{g}_{m})$% , using all points in $P$
% the population $P=\left\{x_{1}, x_{2},\cdots,x_{n}\right\}$
\State Select $\rho_{n} \in \Xi$ and $\epsilon_{i}^{(n)}$ according to COBRA base algorithm
\State $\vec{x}_{start} \leftarrow$ \Call{RandomStart}{$\vec{x}_{best},N_{max}$}		% eliminated parameter c
\State $\vec{x}_{new} \leftarrow$ \Call{OptimCOBRA}{$\vec{x}_{start}, \vec{s}\,^{(n)}$}    	\Comment{see Eq.~\eqref{eq:sub} }
%\Comment{\parbox{5.5cm}{ optimize the surrogate problem of Eq.~\eqref{eq:sub} using COBRA base algorithm}}
%\State Starting from $\vec{x}_{start}$, find $\vec{x}_{new}$ by optimizing the constrained surrogate problem of Eq.~\eqref{eq:sub} with the COBRA base algorithm
\State Evaluate $\vec{x}_{new}$ on the real functions $\widetilde{f},\widetilde{\textbf{g}} $
\If{($|P| \mod 10 == 0$)}															\Comment{every 10th iteration}
	\State $Q \leftarrow$ \Call{AnalysePlogEffect}{$f, P, \vec{x}_{new}$}
\EndIf	
%\If{($\vec{x}_{new}$ infeasible)}
	\State $\vec{x}_{new} \leftarrow$ \Call{repairRI2}{$\vec{x}_{new}$} %, \vec{s}\,^{(n)}(\vec{x}_{new})$}	
								\Comment{\parbox{4.2cm}{  see Koch et al.~\protect\cite{Koch2015a} for details on RI2 (repair algo)}}
%\EndIf
\State $(P,\vec{x}_{best}) \leftarrow $ \Call{updateBest}{$P,\vec{x}_{new},\vec{x}_{best}$}  % eliminated parameter c
\EndWhile
\State \Return $\vec{x}_{best}$
\EndFunction

\item[]
\Function{updateBest}{$P,\vec{x}_{new},\vec{x}_{best}$} % eliminated parameter c
%\State $F_{success} = \textsc{False}$
\State $P \leftarrow P \cup \{\vec{x}_{new}\}$
\If {($\vec{x}_{new} ~\mbox{is feasible} ~~AND~~ \vec{x}_{new} < \vec{x}_{best}$ )}
	%\State $\vec{x}_{best} =  \vec{x}_{new}$
	%\State $F_{success} = \textsc{True}$
	\State \Return ($P,\vec{x}_{new}$)		% eliminated parameter c=0	\Comment{reset counter $c$ if $\vec{x}_{new}$ is better than $\vec{x}_{best}$ }
\EndIf
\State \Return ($P,\vec{x}_{best}$)			% eliminated parameter c+1	\Comment{increment counter $c$ if iteration was unsuccessful }
\EndFunction
\end{algorithmic}
\end{algorithm}%

\begin{algorithm}%
\caption{SACOBRA adjustment functions}
\label{alg:sacobra2Adjust}
	\begin{algorithmic}[1] 
\Function{AnalyseInitialPopulation}{$P,f,\textbf{g}$}
\State $\widehat{FR}=\stackrel[P]{}{\max} {f(P)} -\stackrel[P]{}{\min} {f(P)}$							\Comment{range of objective function}
%\For{$i=1,\ldots,m$} 	
\State $\widehat{GR}_i = \stackrel[P]{}{\max} {g_i(P)} - \stackrel[P]{}{\min} {g_i(P)} \quad\forall i=1,\ldots,m	$  
																																										  			%\Comment{range of each constraint function}
\EndFunction

\item[]
\Function{AdjustConstraintFunction}{$\widehat{GR}_i$,$\textbf{g}$}
%\For{ $i= 1,\ldots m$} 
	\State $g_i() \leftarrow g_i() \cdot \frac{ \mbox{avg}\left( \widehat{GR}_i \right) }{\widehat{GR}_i}  \quad\forall i=1,\ldots,m	$ 									\Comment{see Eq.~\eqref{eq:widehat_GR}} 
%\EndFor
\State \Return $\textbf{g}$
\EndFunction

\item[]
\Function{AdjustDRC}{ $\widehat{FR}$}
\If{$\widehat{FR} > FR_{l}$}					\Comment{Threshold $FR_{l} = 1000$}
\State $\Xi=\Xi_{s}=\left\langle 0.001,0.0\right\rangle$
\Else
\State $\Xi=\Xi_{l}=\left\langle  0.3,0.05, 0.001, 0.0005,0.0 \right\rangle$
\EndIf
\EndFunction

\item[]
\Function{AnalysePlogEffect}{$f, P,\vec{x}_{new}$}			\Comment{$\vec{x}_{new} \notin P$}
\State $S_f \leftarrow$ surrogate model for $f()$ using all points in $P$
\State $S_p \leftarrow$ surrogate model for $plog(f())$ using all points in $P$				
																										\Comment{\parbox{4.2cm}{see function $plog$ in Eq.~\eqref{eq:plog}}}
%\State $error_{S_f} \leftarrow = |S_f(\vec{x}_{new})-f(\vec{x}_{new})|,~ error_{S_p} \leftarrow = |S_p(\vec{x}_{new})-f(\vec{x}_{new})|$
\State $E \leftarrow E\cup \left\{ \dfrac{|S_f(\vec{x}_{new})-f(\vec{x}_{new})|}
																				 {|plog^{-1}\left(S_p(\vec{x}_{new})\right)-f(\vec{x}_{new})|}
													 \right\}	$								\Comment{\parbox{4.5cm}{$E$, the set of approximation error ratios, is initially empty}}
\State \Return $Q=  \log_{10}\left(\med (E) \right)$				% Q is former P_{effect}
\EndFunction

%\item[]
%\Function{AdjustFitnessFunction}{$f,Q$}			% Q is former P_{effect}
%\If{$Q > 1 $}														
%\State \Return $plog(f())$				\Comment{see function $plog$ in Eq.~\eqref{eq:plog}}
%\EndIf
%\State \Return $f()$
%\EndFunction
\end{algorithmic}
\end{algorithm}%

\begin{algorithm}%
\caption{\textsc{RandomStart (RS)}. \textbf{Input}: $\vec{x}_{best}$: the ever-best feasible solution.
%, budget $N_{max}$, $c$: counter for the number of unsuccessful iterations in a row. 
Parameters: restart probabilities $p_1=0.125,~ p_2=0.4$. \textbf{Output}: New starting point $\vec{x}_{start}$. }
\label{alg:rStart2}
	\begin{algorithmic}[1] 
	\Function{RandomStart}{ $\vec{x}_{best}$} %, N_{max},c$}
		%\State $c \leftarrow (F_{success} ~~?~~ 0 ~:~ c + 1 ~)~ $
		\State $F_{low} \leftarrow (|P_{feas}|/|P| < 0.05)$					\Comment{\textsc{True} if less than 5\% of the population are feasible}
		\State $p \leftarrow (F_{low}=TRUE ~~?~~ p_2 ~:~ p_1 ~)~ $
		\State $\epsilon \leftarrow$ a random value $\in [0,1]$ 
		\If{$(\epsilon < p)$ }
		%\If{$(c > N_{max}/10 ~~OR~~ \epsilon < p)$ }
		%\State $c \leftarrow 0$																			\Comment{Reset counter $c$ when choosing a random start point}
		\State $\vec{x}_{start} \leftarrow \mbox{a random point in search space}$
		\Else
		\State $\vec{x}_{start} \leftarrow \vec{x}_{best}$
		\EndIf
		\State \Return $(\vec{x}_{start})$   % eliminated parameter c
\EndFunction

\end{algorithmic}
\end{algorithm}%

%% file: settingTab.tex
\begin{table}[tb]%
\centering
\caption{The default parameter setting used for COBRA. $l$ is the length of the smallest side of the search space (after rescaling, if rescaling is done). The settings for $T_{feas}, T_{infeas}$ proportional to $\sqrt{d}$ ($d$: dimension of problem) are taken from~\cite{regis2014constrained}.
%\SB{I am not sure to keep the last four items or not? If we do not keep them then I think we do not need two columns of parameters because they are exactly same for COBRA and SACOBRA. But I think if we mention plog we should also mention other elements of SACOBRA as parameter}
%\WK{this is o.k.}
}
\label{tab:setting}
\begin{tabular}{lll}
\firsthline
\multirow{2}{*}{parameter} & \multicolumn{2}{c}{value}\\  \cline{2-3}
						& \footnotesize COBRA-R & \footnotesize SACOBRA\\\hline
$\epsilon_{init}$&$0.005 \cdot l$&$0.005 \cdot l$ \\
$\epsilon_{max}$&$ 0.01 \cdot l$&$ 0.01 \cdot l$\\
$T_{feas}$& $  { \lfloor 2\sqrt{d} \rfloor}$&$ { \lfloor 2\sqrt{d} \rfloor}$\\
$T_{infeas}$&$ { \lfloor 2\sqrt{d} \rfloor}$&$ { \lfloor 2\sqrt{d} \rfloor}$\\
$\Xi$&$                 \{ 0.3, 0.05, 0.001, 0.0005, 0.0\}$           & Adaptive\\
$plog(.)$&Never&Adaptive\\
$aCF$&Never&Always\\
$RS$&Never&Adaptive\\

\lasthline
\end{tabular}
\end{table}

%% file: longTab.tex
%\newpage
%\setlength\LTleft{-30pt}            % default: \fill
%\setlength\LTright{-30pt}   
\tiny
%%\singlespacing
%\tablefirsthead{Fct.									& Optimum		  &	& SACOBRA-R            	& COBRA     & ISRES	& RGA 10\%   & COBYLA\\ }
  %\tablehead{Fct.									& Optimum		  &	& SACOBRA-R            	& COBRA     & ISRES	& RGA 10\%   & COBYLA\\ }
  %\tabletail{				& 		  &	&             	&      & 	&   & \\ }
  %\tablelasttail{									& 		  &	&             	&      & 	&   & \\ }
  %\begin{supertabular}{llp{0.3cm}lp{1.8cm}lll}	
	\hspace{-1cm}
	%\rowcolors{1}{}{lightred}

% \begin{longtable}{llp{0.3cm}lp{1.8cm}lll}  
\begin{center}
	
\begin{table}%
\caption{Different optimizers: median (m) of best feasible results and (fe) average number of function evaluations. Results from 30 independent runs with different random number seeds. Numbers in  \textbf{\textcolor[rgb]{0,0,0.55}{boldface (blue)}}: distance to the optimum $\leq 0.001$. Numbers in \textbf{\textcolor[rgb]{1,0,0}{\em italic (red)}}: reportedly better than the true optimum. COBYLA sometimes returns slightly infeasible solutions (number of infeasible runs in brackets).  
%\WK{SB, could you add a column with results for DE?}\SB{Unfortunately, the runs regarding G02 were for 10 dimensional case that's why the regarding cell is empty. the result will be added...}
}
\label{tab:bigTab}
  \resizebox{\textwidth}{!}{  
\begin{tabular}{lr|l|rrrrrr}                      %lllp{2cm}p{2.5cm}}

	\toprule
	\multirow{2}{*}{Fct.}	& \multirow{2}{*}{Optimum}		  &	& SACOBRA         & COBRA      & ISRES 	& RGA 10\%    & COBYLA& DE\\
	& &	& [this work]&\cite{regis2014constrained}& \cite{runarsson2005search,runarsson2000stochastic} & \cite{chootinan2006constraint} & \cite{powell1994direct} (infeas)& \cite{brest2006saDE,Zhang2012}\\
	\toprule
	\multirow{2}{*}{G01}	&	\multirow{2}{*}{-15.0}						
																			&m&\textbf{\textcolor[rgb]{0,0,0.55}{-15.0}}		  &	NA							&\textbf{\textcolor[rgb]{0,0,0.55}{-15.0}} 			&	\textbf{\textcolor[rgb]{0,0,0.55}{-15.0}}		 	& -13.83& \textbf{\textcolor[rgb]{0,0,0.55}{-15.0}}\\
												&							&fe&100			          &	NA							&	350000							& 95512								& 12743&59129\\ %.9\\
	\midrule											
	
	\multirow{2}{*}{G02}	&	\multirow{2}{*}{-0.8036}

												&m&\cellcolor{red!25}	-0.3466					&	\cellcolor{red!25}NA							&	\cellcolor{red!25}-0.7931	&	\cellcolor{red!25}-0.7857			    		& \cellcolor{red!25}-0.197 (5)& \cellcolor{red!25}\textbf{\textcolor[rgb]{0,0,0.55}{-0.8036}}\\
												&							&fe&\cellcolor{red!25}	400				      &\cellcolor{red!25}	NA							&	\cellcolor{red!25}349600					  	&\cellcolor{red!25}	331972								&\cellcolor{red!25}97391& \cellcolor{red!25}226994\\ %.28\\
	\midrule											
	
	\multirow{2}{*}{G03}	&	\multirow{2}{*}{-1.0}

																		&m&\textbf{\textcolor[rgb]{0,0,0.55}{-1.0}}	&	-0.09 					&\textbf{\textcolor[rgb]{1,0,0}{\em -1.001}}		&	-0.9999 	& \textbf{\textcolor[rgb]{0,0,0.55}{-1.0}} (3)&-0.9999\\

												&							&fe&300			            &	100							&	349200				  	&	399804					& 31069&211966\\ %.1\\
  \midrule											
	
	\multirow{2}{*}{G04}	&	\multirow{2}{*}{-30665.539}

																&m&\textbf{\textcolor[rgb]{0,0,0.55}{-30665.539}}&	-30665.15			  &	\textbf{\textcolor[rgb]{0,0,0.55}{-30665.539}}	&	\textbf{\textcolor[rgb]{0,0,0.55}{-30665.539}}							 &\textbf{\textcolor[rgb]{0,0,0.55}{-30665.539}}
																															&\textbf{\textcolor[rgb]{0,0,0.55}{-30665.539}}\\

												&							&fe&200         				&	100							&	192000						&	26981								& 418& 33963\\ %.6\\
	\midrule	

	\multirow{2}{*}{G05}	&	\multirow{2}{*}{5126.497}

																		&m&	\textcolor[rgb]{0,0,0.55}{\textbf{5126.498}}					&	5126.51					&	\textbf{\textcolor[rgb]{0,0,0.55}{5126.497}}	&	\textcolor[rgb]{0,0,0.55}{\textbf{5126.498}}								& \textbf{\textcolor[rgb]{0,0,0.55}{5126.498}} (7)
																															&\textbf{\textcolor[rgb]{0,0,0.55}{5126.498}}\\

												&							&fe&200           			&	100         		&	195600						&	39459								& 194& 13375\\ % .34\\
	\midrule

	\multirow{2}{*}{G06}	&	\multirow{2}{*}{-6961.81}

																		&m&	\textcolor[rgb]{0,0,0.55}{\textbf{-6961.81}}					&	-6834.48				&	\textbf{\textcolor[rgb]{0,0,0.55}{-6961.81}}	&	\textbf{\textcolor[rgb]{0,0,0.55}{-6961.81}}				 & \textbf{\textcolor[rgb]{0,0,0.55}{-6961.81}} (3)
																															&\textbf{\textcolor[rgb]{0,0,0.55}{-6961.81}}\\
	
												&							&fe&100		            	&	100							&	168800						&	13577									& 134& 2857\\
	\midrule											
	
	\multirow{2}{*}{G07}	&	\multirow{2}{*}{24.306}

											  							&m&	\textbf{\textcolor[rgb]{0,0,0.55}{24.306}}&25.32					  &\textbf{\textcolor[rgb]{0,0,0.55}{24.306}}		  &24.471   									& \textbf{\textcolor[rgb]{0,0,0.55}{24.306}} (6)& \textbf{\textcolor[rgb]{0,0,0.55}{24.306}}\\

											  &							&fe&200		            &	100						&	350000						  &	428314								&13072& 94313\\ %.1\\
	\midrule	

	\multirow{2}{*}{G08}	&	\multirow{2}{*}{-0.0958}
		
											  						&m&	\textbf{\textcolor[rgb]{0,0,0.55}{-0.0958}}	&\textbf{\textcolor[rgb]{1,0,0}{\em -0.1}}	&\textbf{\textcolor[rgb]{0,0,0.55}{-0.0958}}					&\textbf{\textcolor[rgb]{0,0,0.55}{-0.0958}}									& -0.0282
																															&\textbf{\textcolor[rgb]{0,0,0.55}{-0.0958}}\\

												&							&fe&200         	&	100							           &	160000						            &	6217									&553&990\\ %.6\\
	\midrule

	\multirow{2}{*}{G09}	&	\multirow{2}{*}{680.630}

											  						&m&680.761			  &3953.97					&\textbf{\textcolor[rgb]{0,0,0.55}{680.630}}			&680.638									    & \textbf{\textcolor[rgb]{0,0,0.55}{680.630}} (2)& \textbf{\textcolor[rgb]{0,0,0.55}{680.630}}\\

												&							&fe&	300						&	100							&	271200						  &	388453								&  8973& 34836\\ %.4\\
	\midrule

	\multirow{2}{*}{G10}	&	\multirow{2}{*}{7049.248}

											  							&m&	7049.253					&18031.74				&\textbf{\textcolor[rgb]{0,0,0.55}{7049.248}}	&7049.566									    & 7064.8 (22)&\textbf{\textcolor[rgb]{0,0,0.55}{7049.248}}\\

												&							&fe&300	          	&	100							&	348800					  &	572629									& 270840& 74875\\ %.2\\
	\midrule

	\multirow{2}{*}{G11}	&	\multirow{2}{*}{0.75}

											  						&m&	\textbf{\textcolor[rgb]{0,0,0.55}{0.75}}		&NA							&\textbf{\textcolor[rgb]{0,0,0.55}{0.75}}					&\textbf{\textcolor[rgb]{0,0,0.55}{0.75}} 							& \textbf{\textcolor[rgb]{0,0,0.55}{0.75}}	&\textbf{\textcolor[rgb]{0,0,0.55}{0.75}}\\

												&							&fe&100           	&NA							&	137200						    &	7215						& 11788& 2190\\ %.7\\
												
	\toprule
					\multicolumn{3}{r|}{average fe}&218	   & 100&  261127	& 210012	& 40652	& 68681 \\
	\midrule
					\multicolumn{3}{r|}{  total fe}&2400	 & 800&	2872400	&2310133	&447175	&755488 \\
								
												%\midrule
%	\multirow{2}{*}{MOPTA}	&	\multirow{2}{*}{$\approx$222}
%
%											  							&m&	222					      &222				&NA	&NA									    & NA\\
%
%												&							&fe&1000	          	&	1000							&	--					  &	--									& --\\
	\bottomrule
	%\end{tabular}}
	\normalsize
	%}	
%\end{longtable}
\end{tabular}
}
\end{table}
\end{center}

\normalsize 
%\singlespace